\documentclass[10pt]{amsart}

\usepackage[margin=1.25in]{geometry} 
\usepackage{amsmath,amsthm,amssymb,bbm,bm}
\usepackage{graphicx}
\usepackage{hyperref}
\usepackage{amssymb,amsfonts}
\usepackage{enumerate, enumitem}
\usepackage{mathrsfs}
\usepackage{xcolor}

\newcommand{\N}{\mathbb{N}}

\newcommand{\R}{\mathbb{R}}

\newcommand{\Fc}{\mathcal{F}}
\newcommand{\E}{\mathbb{E}}

\newcommand{\I}{\mathbbm{1}}

\newcommand{\Ber}{\mathtt{Ber}}
\newcommand{\Beta}{\mathtt{Beta}}
\newcommand{\Poi}{\mathtt{Poisson}}

\newcommand{\PD}{\mathtt{PD}}
\newcommand{\PYP}{\mathtt{PYP}}
\newcommand{\HPYP}{\mathtt{HPYP}}

\newcommand{\EPM}{\mathtt{EPM}}

\newcommand{\SSRM}{\mathtt{SSRM}}
\newcommand{\SSM}{\mathtt{SSM}}
\newcommand{\HSSRM}{\mathtt{HSSRM}}

\newcommand{\Address}{{
\bigskip
\footnotesize

\textsc{Department of Mathematics \& Statistics, McMaster University, Hamilton, ON, L8S 4K1, Canada}\par\nopagebreak
\textit{E-mail address}: \texttt{\{shuifeng, paguyoj\}@mcmaster.ca}
}}

\def\bal#1\eal{\begin{align*}#1\end{align*}}

\newtheorem{theorem}{Theorem}[section]
\newtheorem{lemma}[theorem]{Lemma}
\newtheorem{proposition}[theorem]{Proposition}

\newtheorem{remark}[theorem]{Remark}

\title[Clusters in hierarchical species sampling models]{Asymptotic behavior of clusters in hierarchical species sampling models}
\author{Stefano Favaro, Shui Feng, and J. E. Paguyo}
\date{}


\begin{document}

\begin{abstract}
Consider a random sample of size $N$ from a hierarchical species sampling model. In this paper, we study the large $N$ asymptotic behavior of the number ${\bf K}_N$ of clusters in the random sample and the number ${\bf \widetilde M}_{\ell,N}$ of clusters represented by exactly $\ell$ latent first-level clusters in the first level of the hierarchical model. In particular, we establish almost sure and $L^p$ convergence for ${\bf \widetilde M}_{\ell,N}$, Gaussian fluctuations and the law of the iterated logarithm for ${\bf K}_N$, and large deviation principles for both ${\bf K}_N$ and ${\bf \widetilde M}_{\ell,N}$. Our approach relies on a random sample size (or random index) representation of the number of clusters through the corresponding non-hierarchical species sampling model.
\end{abstract}

\maketitle

\tableofcontents


\section{Introduction} \label{Sec:Introduction}

\subsection{Species Sampling Models}

Let $(\Omega, \Fc, P)$ be the reference probability space, and let $S$ be Polish space, with $M_1(S)$ denoting the space of distributions on $S$ equipped with the weak topology. Consider a generic population of individuals with type (probability) distribution given by $\bm{p} \in \Delta$, where
\begin{displaymath}
\Delta := \left\{ \bm{p} = (p_1, p_2, \ldots) : 0 \leq p_i \leq 1 \text{ and } \sum_{i=1}^\infty p_i = 1 \right\}
\end{displaymath}
is the infinite dimensional simplex. If ${\bf P}=\{P_{k}\}$ is a $\Delta$-valued random variable, and $\{Y_k\}$ is a sequence of iid $S$-valued random variables with diffuse (non-atomic) base distribution $\nu$ in $M_1(S)$, such that ${\bf P}$ is independent of $\{Y_k\}$, then the discrete random probability measure
\bal
\Xi= \sum_{k=1}^\infty P_k \delta_{Y_k}.
\eal
is referred to as a {\em species sampling random measure} \cite{Pitman96}. Since the order of ${\bf P}=\{P_{k}\}$ is is not relevant in $\Xi$, we may assume the $P_{k}$'s are in descending order. See \cite[Chapter 4]{Pit06} for details.

For $n\geq1$, let $(X_1, X_2, \ldots, X_n)$ be a random sample from $\Xi$. Because of the discreteness of $\Xi$, there are ties among the elements of the random sample, meaning that $(X_1, X_2, \ldots, X_n)$ induces a random partition $\Pi_n$ of $[n]:=\{1, 2, \ldots, n\}$ whose {\em blocks}, or {\em clusters}, are the classes induced by the equivalence relation $i\sim j\iff X_{i}=X_{j}$ almost surely. We denote by $|\Pi_n|=k$ the number of clusters of $\Pi_n=\pi_{n}$, and by $|\Pi_{i,n}|=|\pi_{i,n}|$ the number of elements in cluster $i$, for $i=1,\ldots,k$. Following \cite{Pitman96}, the sequence of random partitions $(\Pi_n)_{n \geq 1}$ defines an exchangeable random partition of $\N$, that is, $\Pi_{m}$ is the restriction of $\Pi_{n}$ to the set of partitions of $[m]$, for all $m<n$ and all $n\geq1$, and the distribution of $(\Pi_n)_{n \geq 1}$ is invariant under all finite permutations of $\N$, namely for all $n\geq1$
\begin{displaymath}
P(\Pi_n = \pi_n) = q(|\pi_{1,n}|, \ldots, |\pi_{k,n}|)
\end{displaymath}
for some symmetric function $q$ of $\{|\pi_{1,n}|, \ldots, |\pi_{k,n}|\}$, which is referred to as the {\em exchangeable partition probability function} (EPPF). The distribution of the random partition $\Pi_{n}$ is called a {\em species sampling model}, denoted by $\SSM(q)$. See \cite[Chapter 2]{Pit06} and references therein.

By Kingman's correspondence \cite{Kin78}, the distribution of the species sampling random measure $\Xi$ is determined by the EPPF $q$ and the base distribution $\nu$; see \cite[Chapter 2]{Pit06}. We denote the distribution of $\Xi$ by $\SSRM(q, \nu)$, where $q$ is the EPPF of the random partition induced by random sampling $\Xi$, and $\nu$ is the base distribution of $\Xi$. As a shorthand for the procedure of random sampling from $\Xi$, we say that $(X_1, X_2, \ldots, X_n)$ is a random sample of size $n$ from $\SSRM(q, \nu)$, and that the induced random partition of $[n]$ is distributed as $\SSM(q)$.

\subsubsection{Ewens-Pitman Sampling Model}

The {\em Pitman-Yor process} \cite{Per(92),Pit(95),PY97} is the most popular species sampling random measure. For a pair of parameters $\alpha \in [0,1)$ and $\theta > -\alpha$, and a base distribution $\nu$, the Pitman-Yor process (PYP), whose distribution is denoted by $\PYP(\alpha, \theta, \nu)$, is the discrete random measure
\bal
\Xi_{\alpha, \theta, \nu} = \sum_{k=1}^\infty V_k \delta_{Y_k}
\eal
where $\bm{V} = (V_1, V_2, \ldots)$, rearranged in weakly decreasing order, is distributed according to the {\em two-parameter Poisson-Dirichlet distribution}, denoted by $\PD(\alpha, \theta)$. When $\alpha = 0$, the PYP reduces to the {\em Dirichlet process} \cite{Fer73}. See \cite{Fen10} and references therein for a comprehensive account on the Poisson-Dirichlet distribution and its two-parameter generalization.

Let $(X_1, \ldots, X_n)$ be a random sample of size $n$ from $\PYP(\alpha, \theta, \nu)$. The distribution of the induced random partition is the {\em Ewens-Pitman sampling model} with parameters $\alpha \in [0,1)$ and $\theta > -\alpha$, denoted by $\EPM(\alpha, \theta)$ \cite{Pit(95)}. When $\alpha = 0$, one recovers the {\em Ewens model} \cite{Ewe(72)}. We refer to Pitman \cite[Chapter 3]{Pit06} and Crane \cite{Cra16} for an overview of the Ewens-Pitman sampling model. 

Let $\Pi_n$ be a random partition of $[n]$ distributed as $\EPM(\alpha, \theta)$. Let 
\begin{displaymath}
K_n(\alpha,\theta) := |\Pi_n|
\end{displaymath}
be the {\em number of clusters} in $\Pi_n$, and let 
\begin{displaymath}
M_{r,n}(\alpha, \theta) := \sum_{i=1}^{K_n(\alpha,\theta)} \I_{\{|\Pi_{i,n}| = r\}}, \qquad \text{for $1\leq r\leq n$,}
\end{displaymath}
be the {\em number of clusters of size $r$} in $\Pi_{n}$, where $\I_A$ is the indicator random variable of $A$. Equivalently, $K_n(\alpha,\theta)$ is the {\em number of distinct values} and $M_{r,n}(\alpha,\theta)$ is the {\em number of distinct values with frequency $r$}, for $1 \leq r \leq n$, in a random sample $(X_1, \ldots, X_n)$ of size $n$ from $\PYP(\alpha, \theta, \nu)$. If $\alpha = 0$, $K_n(\theta) := K_n(0, \theta)$ and $M_{r,n}(\theta):=M_{r,n}(0, \theta)$. 

There is an extensive literature investigating the large $n$ asymptotic behavior of $K_n(\alpha,\theta)$ and $M_{r,n}(\alpha, \theta)$, for $\alpha\in[0,1)$ and $\theta>-\alpha$. For $\alpha=0$ and $\theta>0$, \cite[Theorem 2.3]{KH73} showed that 
\begin{equation}\label{old-1} 
\lim_{n \rightarrow \infty}\frac{K_n(\theta)}{\log n} = \theta
\end{equation}
almost surely. Moreover, since $K_n(\theta)$ can be represented as the sum of independent Bernoulli random variables \cite[Chapter 4]{Pit06}, by means of the Lindeberg-Feller central limit theorem 
 \begin{equation}\label{old-2} 
\frac{K_n(\theta) - \theta \log n}{\sqrt{\theta \log n}} \rightarrow N(0,1),\qquad  n\rightarrow \infty,
\end{equation}
where the convergence is in distribution and $N(0,1)$ is a standard normal random variable. As a generalization of these results to the $M_{r,n}(\alpha, \theta)$'s, \cite[Theorem 1]{ABT92} showed that
\begin{equation}\label{old-3} 
(M_{1,n}(\theta), M_{2,n}(\theta), \ldots)\rightarrow (Z_1, Z_2, \ldots),\qquad  n\rightarrow \infty,
\end{equation}
where the convergence is in distribution and $\{Z_r\}$ are independent Poisson random variables such that $Z_r$ is $\Poi(\theta/r)$ distributed, for $r\geq1$. See also \cite{Arr(03)} and references therein.

As a natural generalization of \eqref{old-1}, for $\alpha \in (0,1)$ and $\theta > -\alpha$, \cite[Theorem 3.8]{Pit06} showed that 
\begin{equation}\label{old-4} 
\lim_{n\rightarrow \infty}\frac{K_n(\alpha,\theta)}{n^\alpha} = S_{\alpha, \theta},
\end{equation}
almost surely and in $L^p$, where $S_{\alpha, \theta}$ is a scaled Mittag-Leffler random variable commonly referred to as {\em Pitman's $\alpha$-diversity}; see \cite[Chapter 3]{Pit06}. Furthermore, from \cite[Lemma 3.11]{Pit06},
\begin{equation}\label{old-5} 
\lim_{n \rightarrow\infty}\frac{M_{r,n}(\alpha,\theta)}{n^\alpha} \to p_\alpha(r)S_{\alpha, \theta}
\end{equation}
almost surely and in $L^p$, where $p_\alpha(r) := \frac{\alpha(1-\alpha)^{(r-1)}}{r!}$ and $(n)^{(k)} := n(n+1)\dotsb (n+k-1)$ is the rising factorial. We refer to \cite{BF24} for a martingale approach to the derivation of \eqref{old-4} and \eqref{old-5}. A large deviation principle for $K_{n}(\alpha,\theta)$ was first established in \cite{FH98}, and then extended to $M_{r,n}(\alpha,\theta)$ in \cite{FF15}. Moderate deviation principles for $K_n(\alpha,\theta)$ and $M_{r,n}(\alpha,\theta)$ are in \cite{FFG18}. 

Further results on $K_n(\alpha,\theta)$ and $M_{r,n}(\alpha,\theta)$ include Berry-Essen type inequalities \cite{DF20}, Gaussian fluctuations \cite{BF24}, and concentration inequalities \cite{Per(22),BF25}. See also \cite{MS25} for an analysis on the large $n$ asymptotic behavior of the joint distribution of $(K_n(\alpha, \theta), M_{r,n}(\alpha, \theta))$.

\subsection{Hierarchical Species Sampling Models} \label{Sec:HSSM}

Hierarchical species sampling models generalize species sampling models to multiple populations of individuals. Fix a number $d \geq 1$ of populations or groups. A {\em hierarchical species sampling random measure} is a vector of $d$ discrete random probability measures $(\Xi_1, \ldots, \Xi_d)$ such that
\begin{eqnarray}\label{SS-Model}
\Xi_i \mid \Xi_0 &\stackrel{iid}{\sim}& \SSRM(q, \Xi_0), \qquad 1 \leq i \leq d \\
\Xi_0 &\sim& \SSRM(q_0, \nu), \nonumber
\end{eqnarray}
where $\sim$ denotes distributed as, $\nu$ is the base distribution, and $q_0, q$ are two EPPFs.  In particular, $\Xi_0$  is a discrete random probability measure describing the global distribution of the hierarchical model, while the discrete random probability measures $\Xi_1, \ldots, \Xi_d$ describe the local distributions of the $d$ groups. All $d$ groups interact with each other through the global distribution $\Xi_0$, and hence we get the term hierarchical. The distribution of $(\Xi_1, \ldots, \Xi_d)$  is denoted by $\HSSRM(q_0, q, \nu, d)$, where $q_0$ is the global EPPF and $q$ is the local EPPF.

We say that {\em $(X_{i,j})_{1 \leq i \leq d, j \geq 1}$ is a random sample from $\HSSRM(q_0, q, \nu, d)$} if its elements are conditionally independent given $(\Xi_1, \ldots, \Xi_d)$ such that $X_{i,j} \mid (\Xi_1, \ldots, \Xi_d) \stackrel{iid}{\sim} \Xi_i$, for all $1\leq i \leq d$ and $j \geq 1$. In particular, by means of de Finetti's theorem, $(X_{i,j})_{1 \leq i \leq d, j \geq 1}$ consists of $d$ {\em partially exchangeable} random samples, which means random samples within a fixed row are {\em exchangeable}, but not across different rows. That is, $\{(X_{i,j})_{j \geq 1}\} = \{(X_{i,\sigma_i(j)})_{j \geq 1}\}$ in distribution, for all finite permutations $(\sigma_1, \ldots, \sigma_d)$. For hierarchical species sampling models, we use capital letters $N, N_i$ to denote the corresponding sample sizes. 

Consider an array $(X_{i,j})_{1 \leq i \leq d, 1 \leq j \leq N_i}$ of size $N = \sum_{i=1}^d N_i$ consisting of $d$ partially exchangeable sequences randomly sampled from $\HSSRM(q_0, q, \nu, d)$. More specifically, a random sample of size $N$ from $\HSSRM(q_0, q, \nu, d)$ consists of samples from each of the $d$ groups. The sample size from group $i$ is $N_i = N_i(N)$, which depends on $N$. Since the random measures $(\Xi_1, \ldots, \Xi_d)$ and $\Xi_0$ are discrete, $P(X_{\ell,i} = X_{k,j}) > 0$ for any $\ell$ and $k$, so that there is a positive probability of equality, both within each sample and across the different samples. Accordingly, the array $(X_{i,j})_{1 \leq i \leq d, 1 \leq j \leq N_i}$ induces a random partition of $[N]$, such that, for any $1\leq \ell, k \leq d, 1\leq \jmath \leq N_\ell, 1\leq j\leq N_k$, $\left( \sum_{m=1}^{\ell - 1} N_m \right) + \jmath$ and $\left( \sum_{m=1}^{k - 1} N_m \right) + j$ are in the same cluster if and only if $X_{\ell,\jmath} = X_{k,j}$ almost surely. The distribution of the induced random partition is called a {\em hierarchical species sampling model}.

\subsubsection{Hierarchical Ewens-Pitman Sampling Model}

The main focus of this paper is on the {\em hierarchical Pitman-Yor process}, which is a hierarchical species sampling random measure \cite{TJ10}. More specifically, for a fixed $d \geq 1$, the hierarchical Pitman-Yor process is a vector of $d$ random measures $(\Xi_1, \ldots,\Xi_d)$ satisfying   
 \begin{eqnarray}\label{PY-Model}
\Xi_i \mid \Xi_0 &\stackrel{iid}{\sim} &\PYP(\beta, \theta, \Xi_0), \qquad i \in [d] \\
\Xi_0 &\sim& \PYP(\alpha, \theta_0, \nu),\nonumber
\end{eqnarray}
where $\alpha, \beta \in [0,1)$, $\theta > -\beta$, and $\theta_0 > -\alpha$. By setting $\alpha = \beta = 0$, one recovers the {\em hierarchical Dirichlet process}, first introduced in \cite{TJBB06}. We can also get mixed cases by letting one level of the hierarchy be a Dirichlet process and the other a Pitman-Yor process. We denote by $\HPYP(\alpha, \theta_0, \beta, \theta, \nu, d)$ the distribution of the hierarchical Pitman-Yor process. 

In analogy to $\PYP(\alpha, \theta, \nu)$, let $(X_{i,j})_{1 \leq i \leq d, 1 \leq j \leq N_i}$ be a random sample of size $N = \sum_{i=1}^d N_i$ from $\HPYP(\alpha, \theta_0, \beta, \theta, \nu, d)$. The distribution of the induced random partition is called the {\em hierarchical Ewens-Pitman sampling model}. See \cite{TJBB06,TJ10,CLOP19} and references therein.

\subsection{Main Results} \label{Sec:MainResults}

To introduce our main results, for $d \geq 1$ let $(X_{i,j})_{1 \leq i \leq d, 1 \leq j \leq N_i}$ be a random sample of size $N = \sum_{i=1}^d N_i$ from $\HPYP(\alpha, \theta_0, \beta, \theta, \nu, d)$. There are several important statistics associated with the random sample, namely:
\begin{itemize}
\item ${\bf K}_N(\alpha,\theta_0; \beta,\theta)$, the {\em number of distinct values in the random sample}, and 
\item ${\bf M}_{r,N}(\alpha,\theta_0; \beta,\theta)$, the {\em number of distinct values with frequency $r$ in the random sample}, for any $1 \leq r \leq N$. 
\item ${\bf \widetilde M}_{\ell,N}(\alpha,\theta_0;\beta,\theta)$, the {\em number of distinct values represented by exactly $\ell$ clusters in the first level of the hierarchical Pitman-Yor process}, for any $1 \leq \ell \leq N$, which we define rigorously in Section \ref{sec_randind}. 
\end{itemize}
Let $\Pi_N = (\Pi_{1,N}, \ldots, \Pi_{{\bf K}_N(\alpha,\theta_0; \beta,\theta), N})$ be the random partition of $[N]$ induced by the random sample. Then we can also define ${\bf K}_N(\alpha,\theta_0; \beta,\theta)$ and ${\bf M}_{r,N}(\alpha,\theta_0; \beta,\theta)$ in terms of clusters in $\Pi_N$: 
\begin{displaymath}
{\bf K}_N(\alpha,\theta_0; \beta,\theta) = |\Pi_N|
\end{displaymath}
is the number of clusters in $\Pi_N$, and 
\begin{displaymath}
{\bf M}_{r,N}(\alpha,\theta_0; \beta,\theta) = \sum_{k=1}^{{\bf K}_N(\alpha,\theta_0; \beta,\theta)} \I_{\{|\Pi_{k,N}| = r\}}, \qquad 1\leq r\leq N,
\end{displaymath}
is the number of clusters of size $r$ in $\Pi_N$. The statistic ${\bf \widetilde M}_{\ell,N}(\alpha,\theta_0;\beta,\theta)$ is not a function of the observed random sample alone, but also depends on the latent first level clustering induced by the hierarchical Pitman-Yor process, and so we postpone its formal definition until Section \ref{sec_randind} (see equation (\ref{tablefrequency})). 

The statistics  ${\bf M}_{r,N}(\alpha,\theta_0; \beta,\theta)$ and ${\bf \widetilde M}_{\ell,N}(\alpha,\theta_0;\beta,\theta)$ describe different partition structures in the hierarchical model. For example, consider the second level labels as surnames of individuals. Given $N$ individuals,  ${\bf M}_{r,N}(\alpha,\theta_0; \beta,\theta)$ represents the number of distinct surnames, each of which correspond to exactly $r$ individuals. Consider each cluster in the first level of the hierarchy as a family in which all members have the same surname. Then ${\bf \widetilde M}_{\ell,N}(\alpha,\theta_0;\beta,\theta)$ represents the number of distinct surnames each of which is shared by exactly $\ell$ families.

There are four special cases, depending on whether or not $\alpha = 0$ or $\beta = 0$. We use the following notation to distinguish between these cases:
\begin{itemize}
\item[(i)] If $\alpha = \beta = 0$, then ${\bf K}_N(\theta_0;\theta):={\bf K}_N(0,\theta_0; 0,\theta)$.
\item[(ii)] If $\alpha = 0$ and $\beta \in (0,1)$, then ${\bf K}_N(\theta_0; \beta,\theta):={\bf K}_N(0,\theta_0; \beta,\theta)$.
\item[(iii)] If $\alpha \in (0,1)$ and $\beta = 0$, then ${\bf K}_N(\alpha, \theta_0; \theta):={\bf K}_N(\alpha,\theta_0; 0,\theta)$
\item[(iv)] Otherwise we use the general notation ${\bf K}_N(\alpha,\theta_0; \beta,\theta)$.
\end{itemize}
The corresponding notation for ${\bf M}_{r,N}(\alpha,\theta_0; \beta,\theta)$ and ${\bf \widetilde M}_{\ell,N}(\alpha,\theta_0; \beta,\theta)$ are defined analogously. 

There has been some recent work on hierarchical species sampling models. A distribution theory for hierarchical random measures has been developed in \cite{CLOP19}, which encompasses both the hierarchical Pitman-Yor and Dirichlet processes, providing the distribution of the number of clusters. Subsequently, \cite{BCR20} obtained finite sample and asymptotic distributions of the number of clusters; see also \cite{BL20, BL21}. We contribute to this recent line of work by studying the large $N$ asymptotic behavior of ${\bf K}_N(\alpha,\theta_0; \beta,\theta)$ and ${\bf \widetilde M}_{\ell,N}(\alpha,\theta_0;\beta,\theta)$. In particular, we establish various limit theorems, including almost sure and $L^p$ convergence, Gaussian fluctuations, and large deviation principles. Our main results are:
\begin{itemize}[leftmargin=0.5cm]
\item almost sure and $L^p$ convergence of ${\bf \widetilde M}_{\ell,N}(\alpha,\theta_0;\beta,\theta)$ (Proposition~\ref{HEPblockLLN}); 
\item Poisson limit theorem for ${\bf \widetilde M}_{\ell,N}(0,\theta_0; \beta, \theta)$ (Theorem~\ref{HEPblockPoissonLimit});
\item Gaussian fluctuation theorem for the number of clusters in a random partition of random size $\xi(N)$ distributed as $\SSM(q_0)$ (Theorem~\ref{ClustersCLT});
\item as an application of Theorem~\ref{ClustersCLT}: Gaussian fluctuation theorems for ${\bf K}_N(\alpha, \theta_0; \beta, \theta)$ in all four cases (Theorem \ref{HEPCLT});
\item large deviation principle for the number of clusters in a random partition of random size $\xi(N)$ distributed as $\SSM(q_0)$ (Theorem~\ref{ClustersLDP});
\item as an application of Theorem~\ref{ClustersLDP}: large deviation principles for ${\bf K}_N(\alpha,\theta_0; \beta,\theta)$ and  ${\bf \widetilde M}_{\ell,N}(\alpha,\theta_0; \beta,\theta)$ in all cases (Theorems~\ref{HEPLDP} and \ref{HEPblockLDP});
\item the law of the iterated logarithm for ${\bf K}_N(\alpha,\theta_0; \beta,\theta)$ (Theorem \ref{HEPMLIL}). 
\end{itemize}

Our proof strategy relies on a random sample size (or random index) representation of the number of clusters through the corresponding non-hierarchical species sampling model. Finally, although ${\bf M}_{r,N}(\alpha,\theta_0; \beta,\theta)$ is a natural statistic of interest, its asymptotic behavior presents greater challenges as it depends on both the joint distribution of the latent first level clustering and the observed partition. Thus it cannot be analyzed through the random index representation, and so a separate treatment is required and will be the subject of forthcoming work. 

\subsection{Outline} \label{Sec:Outline}

In Section \ref{Sec:Prelims}, we introduce the notation and auxiliary results that we use throughout the paper, including the random index representation. Section \ref{Sec:ClustersFrequencyLLN} contains the almost sure and $L^p$ convergence for ${\bf \widetilde M}_{\ell,N}(\alpha,\theta_0; \beta,\theta)$. Section \ref{Sec:GaussianFluctuations} includes a Gaussian fluctuation theorem for the number of clusters in a random partition of random size distributed as $\SSM(q_0)$, which is used to obtain the Gaussian fluctuation theorem for ${\bf K}_N(\alpha, \theta_0; \beta, \theta)$. Large deviation principles are obtained in Section \ref{Sec:LDP} for the number of clusters in a random partition of random size distributed as $\SSM(q_0)$, and for both ${\bf K}_N(\alpha,\theta_0; \beta,\theta)$ and ${\bf \widetilde M}_{\ell,N}(\alpha,\theta_0; \beta,\theta)$. The law of the iterated logarithm for ${\bf K}_N(\alpha,\theta_0; \beta,\theta)$ is established in Section \ref{Sec:LIL}. We conclude the paper by discussing some directions for future work in Section \ref{Sec:FinalRemarks}.


\section{Preliminaries} \label{Sec:Prelims}

\subsection{Notation and Definitions} \label{Sec:Notation}

Let $a_n, b_n$ be two sequences: i) if $\lim_{n \to \infty} \frac{a_n}{b_n}$ is a nonzero constant, then we write $a_n \asymp b_n$ and we say that $a_n$ is of the same {\em order} as $b_n$; ii) if  the limit is $1$, we write $a_n\sim b_n$ and we say that $a_n$ is {\em asymptotic} to $b_n$; iii) if there exists positive constants $c$ and $n_0$ such that $a_n \leq cb_n$ for all $n \geq n_0$, then we write $a_n = O(b_n)$ and say that $a_n$ is {\em big-O} of $b_n$; (iv) if $\lim_{n \to \infty} \frac{a_n}{b_n} = 0$, then we write $a_n = o(b_n)$ and we say that $a_n$ is {\em little-O} of $b_n$. 

We let {\em iid} abbreviate independent and identically distributed. Convergence in distribution, convergence in probability, and convergence almost surely are denoted by $\xrightarrow{D}$, $\xrightarrow{P}$, and $\xrightarrow{a.s.}$, respectively. Similarly equality in distribution, in probability, and almost surely are denoted by $\stackrel{D}{=}$, $\stackrel{P}{=}$, and $\stackrel{a.s.}{=}$, respectively. If $X$ is a random variable distributed as $\nu$, then we write $X \sim \nu$. 

Let $x \in \R$ and $n \in \N$. The {\em rising factorial} of $x$ of order $n$ is denoted by $(x)^{(n)} := x(x+1)\dotsb (x+n-1)$.  

A function $f: (0,\infty) \to (0, \infty)$ is {\em regularly varying} if it is of the form $f(x) = x^\sigma L(x)$, where $\sigma \in \R$ is the {\em index of regular variation} and $L$ is a {\em slowly varying function}, which is a function $L:(0, \infty) \to (0,\infty)$ such for all $a > 0$, $\frac{L(ax)}{L(x)} \to 1$ as $x \to \infty$ \cite{Bin(89)}. For example, $\log n$ is slowly varying and $n^\sigma \log n$ is regularly varying.

Further terminology and notation will be introduced when necessary.

\subsection{Random Index Representation}\label{sec_randind}

For $d\geq1$, consider a random sample $(X_{i,j})_{1 \leq i \leq d, 1 \leq j \leq N_i}$ of size $N = \sum_{i=1}^d N_i$ from $\HPYP(\alpha, \theta_0, \beta, \theta, \nu, d)$. The random sample admits a Chinese Restaurant Franchise representation \cite{TJBB06, CLOP19, BCR20}. In the Chinese Restaurant Franchise representation, the $N$ observations correspond to ``customers" and the $d$ groups correspond to ``restaurants" which comprise the ``franchise". In each restaurant, customers are clustered at the first level of the hierarchy according to ``tables". The tables are then clustered across all restaurants at the second level of the hierarchy according to ``dishes". In particular, the first customer who sits at a table picks a dish from a common franchise-wide menu, and this dish is shared by all other customers who join the same table afterwards. Accordingly, $X_{i,j}$ represents the dish served to the $j$th customer in the $i$th restaurant. 

Following \cite{CLOP19}, we can construct the hierarchical Pitman-Yor process on an enlarged probability space as follows. 
For each $1 \leq i \leq d$, let 
\bal
T_{i,1}, T_{i,2}, \ldots \mid G_i \stackrel{iid}{\sim} G_i,
\eal
where $G_i \stackrel{iid}{\sim} \PYP(\beta, \theta, \tilde{\nu})$ and $\tilde{\nu}$ is a diffuse base measure. Independently of $(G_i)_{1 \leq i \leq d}$, let 
\bal
Z_1, Z_2, \ldots \mid G_0 \stackrel{iid}{\sim} G_0
\eal
where $G_0 \sim \PYP(\alpha, \theta_0, \nu)$. 
For each $1 \leq i \leq d$, let $K_{i,N_i}(\beta,\theta)$ be the number of distinct values in $(T_{i,j})_{1 \leq j \leq N_i}$. Let
\begin{equation}\label{totalclusters}
\xi(N) = \sum_{i=1}^d K_{i,N_i}(\beta,\theta)
\end{equation}
be the total number of distinct values in $(T_{i,j})_{1 \leq i \leq d, 1 \leq j \leq N_i}$. Finally, note that $K_n(\alpha, \theta_0)$ is the number of distinct values in $\{Z_1, \ldots, Z_n\}$.

In the Chinese Restaurant Franchise representation, $T_{i,j}$ is the label of the table where the $j$th customer of the $i$th restaurant is seated, $K_{i,N_i}(\beta,\theta)$ is the number of tables where the $N_i$ customers of restaurant $i$ are seated, $\xi(N)$ is the total number of tables in the entire franchise, and $K_{\xi(N)}(\alpha,\theta_0)$ is the number of distinct dishes served amongst the $\xi(N)$ tables where the $N$ customers in the entire franchise are seated. 

Let $\Pi_N = (\Pi_{1,N}, \ldots, \Pi_{{\bf K}_N(\alpha,\theta_0; \beta,\theta), N})$ be the random partition of $[N]$ induced by a random sample $(X_{i,j})_{1 \leq i \leq d, 1 \leq j \leq N_i}$ of size $N = \sum_{i=1}^d N_i$ from $\HPYP(\alpha, \theta_0, \beta, \theta, \nu, d)$. We refer to the partition sets $(\Pi_1, \ldots, \Pi_{{\bf K}_N(\alpha,\theta_0; \beta,\theta)})$ as the {\em second-level clusters} of $\Pi_N$. Let $\Gamma_N = (\Gamma_{1,N}, \ldots, \Gamma_{\xi(N), N})$ be the latent random partition of $[N]$ induced by the same random sample such that, for any $1\leq \ell, k \leq d, 1\leq \jmath \leq N_\ell, 1\leq j\leq N_k$, $\left( \sum_{m=1}^{\ell - 1} N_m \right) + \jmath$ and $\left( \sum_{m=1}^{k - 1} N_m \right) + j$ are in the same cluster if and only if $T_{\ell,\jmath} = T_{k,j}$ almost surely. We refer to the partition sets $(\Gamma_{1,N}, \ldots, \Gamma_{\xi(N), N})$ as the {\em first-level clusters} of $\Pi_N$.
For each $1 \leq k \leq {\bf K}_N(\alpha,\theta_0; \beta,\theta)$, let $L_{k,N} = |\{j \in [\xi(N)] : \Gamma_{j,N} \subseteq \Pi_{k,N}\}|$ be the number of first-level clusters that belong to the $k$th second-level cluster. Then
\begin{equation} \label{tablefrequency}
{\bf \widetilde M}_{\ell,N}(\alpha, \theta_0; \beta, \theta) = \sum_{k=1}^{{\bf K}_N(\alpha,\theta_0; \beta,\theta)} \I_{\{L_{k,N} = \ell\}}, \qquad 1 \leq \ell \leq N,
\end{equation}
is the {\em number of distinct values represented by exactly $\ell$ first-level clusters}. 

The next proposition establishes the random-index representations for ${\bf K}_N(\alpha,\theta_0; \beta,\theta)$ and ${\bf \widetilde M}_{\ell,N}(\alpha, \theta_0; \beta, \theta)$ (see also \cite[Theorem 6 and Remark 1]{CLOP19}).

\begin{proposition} \label{PropRandIndex}
On an enlarged probability space supporting the construction above, let $(T_{i,j})_{1 \leq i \leq d, 1 \leq j \leq N_i}$ and $(Z_k)_{k \geq 1}$ be the latent random variables as defined above. Let $\tau_1, \ldots, \tau_{\xi(N)}$ denote the distinct values in $(T_{i,j})_{1 \leq i \leq d, 1 \leq j \leq N_i}$, where $\xi(N)$ is the number of distinct values. Define $(R_{i,j})_{1 \leq i \leq d, 1 \leq j \leq N_i}$ such that $R_{i,j} = k$ if and only if $T_{i,j} = \tau_k$. Finally, define the observations $(X_{i,j})_{1 \leq i \leq d, 1 \leq j \leq N_i}$ such that $X_{i,j} = Z_{R_{i,j}}$. 
Then $(X_{i,j})_{1 \leq i \leq d, 1 \leq j \leq N_i}$ is a random sample of size $N = \sum_{i=1}^d N_i$ from $\HPYP(\alpha, \theta_0, \beta, \theta, \nu, d)$, and moreover, under this coupling, 
\begin{equation}\label{RandomIndex-1}
{\bf K}_N(\alpha,\theta_0; \beta,\theta) \stackrel{a.s.}{=} K_{\xi(N)}(\alpha,\theta_0).
\end{equation}
and
\begin{equation}\label{RandomIndex-2}
{\bf \widetilde M}_{\ell,N}(\alpha, \theta_0; \beta, \theta) \stackrel{a.s.}{=} M_{\ell,\xi(N)}(\alpha,\theta_0).
\end{equation}
\end{proposition}

\begin{proof}
The fact that $(X_{i,j})_{1 \leq i \leq d, 1 \leq j \leq N_i}$ is a random sample of size $N = \sum_{i=1}^d N_i$ from $\HPYP(\alpha, \theta_0, \beta, \theta, \nu, d)$ follows immediately from the Chinese Restaurant Franchise representation. Next, since $\tilde{\nu}$ is diffuse, the distinct values in $(\tau_1, \ldots, \tau_{\xi(N)})$ are in bijection with the first-level clusters induced by $(T_{i,j})_{1 \leq i \leq d, 1 \leq j \leq N_i}$. Thus $R_{i,j}$ is the index of the first-level cluster containing the $(i,j)$th observation. By the coupling $X_{i,j} = Z_{R_{i,j}}$, every observation in the same first-level cluster has the same value since $R_{i,j} = R_{k,\ell}$ implies that $X_{i,j} = Z_{R_{i,j}} = Z_{R_{k,\ell}} = X_{k,\ell}$. Moreover, two first-level clusters, say indexed by $r$ and $s$, have the same value if and only if $Z_r = Z_s$. 

Therefore, the distinct values in $(X_{i,j})_{1 \leq i \leq d, 1 \leq j \leq N_i}$ are exactly the distinct values in $(Z_1, \ldots, Z_{\xi(N)})$, from which it follows that
\bal
{\bf K}_N(\alpha,\theta_0; \beta,\theta) \stackrel{a.s.}{=} K_{\xi(N)}(\alpha,\theta_0).
\eal
Moreover, a distinct value in $(X_{i,j})_{1 \leq i \leq d, 1 \leq j \leq N_i}$ is represented by exactly $\ell$ first-level clusters if and only if the corresponding value occurs exactly $\ell$ times in $(Z_1, \ldots, Z_{\xi(N)})$, so that
\bal
&{\bf \widetilde M}_{\ell,N}(\alpha, \theta_0; \beta, \theta) \stackrel{a.s.}{=} M_{\ell,\xi(N)}(\alpha,\theta_0). \qedhere
\eal
\end{proof}

We refer to (\ref{RandomIndex-1}) and (\ref{RandomIndex-2}) as the {\em random index representations} of ${\bf K}_N(\alpha,\theta_0; \beta,\theta)$ and ${\bf \widetilde M}_{\ell,N}(\alpha, \theta_0; \beta, \theta)$, respectively. These play key roles in establishing the asymptotic results in the sequel.

\subsection{Auxiliary Results} \label{Sec:AuxiliaryResults}

In this section, we record six auxiliary results that will be used later in the paper to establish our main results. The first result is the almost sure and $L^p$ convergence for ${\bf K}_{N}(\alpha,\theta_{0};\beta,\theta)$ \cite{BCR20}.

\begin{proposition}\cite[Proposition 9]{BCR20} \label{HEPLLN}
Let ${\bf K}_N(\alpha, \theta_0; \beta, \theta)$ be the number of distinct values in a random sample $(X_{i,j})_{1 \leq i \leq d, 1 \leq j \leq N_i}$ of size $N = \sum_{i=1}^d N_i$ from $\HPYP(\alpha, \theta_0, \beta, \theta, \nu, d)$.
Suppose that $\frac{N_i}{N} = \frac{N_i(N)}{N} \to w_i > 0$ as $N \to \infty$ for all $i \in [d]$. Then as $N \to \infty$, the following almost sure convergences hold:
\begin{enumerate}
\item If $\alpha = \beta = 0$, 
\bal
\frac{{\bf K}_N(\theta_0; \theta)}{\log \log N} \xrightarrow{a.s.} \theta_0,
\eal
\item If $\alpha = 0$ and $\beta \in (0,1)$,
\bal 
\frac{{\bf K}_N(\theta_0; \beta, \theta)}{\beta \log N} \xrightarrow{a.s.} \theta_0, 
\eal
\item If $\alpha \in (0,1)$ and $\beta = 0$,
\bal
\frac{{\bf K}_N(\alpha, \theta_0; \theta)}{(\log N)^\alpha} \xrightarrow{a.s.} S_{\alpha, \theta_0} (\theta d)^\alpha, 
\eal
\item If $\alpha, \beta \in (0,1)$, 
\bal
\frac{{\bf K}_N(\alpha, \theta_0; \beta, \theta)}{N^{\alpha\beta}} \xrightarrow{a.s.} S_{\alpha, \theta_0} \eta^\alpha,
\eal
\end{enumerate}
where $\eta = \sum_{i=1}^d w_i^\beta S_{\beta,\theta,i}$ such that $\{S_{\beta, \theta, i}\}_{i \in [d]}$ are iid copies of $S_{\beta,\theta}$ independently of $S_{\alpha, \theta_0}$. Moreover the above convergences also hold in $L^p$ for all integers $p \geq 1$. 
\end{proposition}

The second result is a Gaussian fluctuation for $K_n(\alpha, \theta)$ in the Ewens-Pitman sampling model \cite{BF24}.

\begin{proposition}\cite[Theorem 2.3]{BF24}\label{EPCLT}
Let $K_n(\alpha, \theta)$ be the number of clusters in a random partition of $[n]$ distributed as $\EPM(\alpha, \theta)$, with $\alpha \in (0,1)$ and $\theta > -\alpha$. Then
\bal
\sqrt{n^\alpha}\left( \frac{K_n(\alpha, \theta)}{n^\alpha} - S_{\alpha, \theta} \right) \stackrel{D}{\longrightarrow} \sqrt{S_{\alpha, \theta}'} N(0,1), 
\eal
where $S_{\alpha, \theta}'$ is a random variable independent of $N(0,1)$ and equals in distribution to $S_{\alpha, \theta}$.
\end{proposition}

The third result is a large deviation principle for $K_n(\alpha, \theta)$ in the Ewens-Pitman sampling model \cite{FH98}.

\begin{proposition}\cite[Theorem 1.2]{FH98} \label{ClustersLDPLevel0}
For $\alpha \in (0,1)$ and $\alpha > -\theta$, the sequence $\left\{\frac{K_n(\alpha, \theta)}{n} \right\}$ satisfies a large deviation principle with speed $n$ and rate function
\bal
I^\alpha(x) = \sup_{\lambda \in \R} \{ \lambda x - \Lambda_\alpha(\lambda)\},
\eal
where
\bal
\Lambda_\alpha(\lambda) = \begin{cases}
-\log[1 - (1 - e^{-\lambda})^{1/\alpha}] & \text{if $\lambda > 0$},\\
0 & \text{otherwise}.
\end{cases}
\eal

For $\alpha = 0$ and $\theta > 0$, $\left\{\frac{K_n(\theta)}{\log n} \right\}$ satisfies a large deviation principle with speed $\log n$ and rate function
\bal
I_\theta(x) = \begin{cases}
x\log \frac{x}{\theta} - x + \theta & \text{for $x \geq 0$} \\
\infty & \text{otherwise}.
\end{cases}
\eal
\end{proposition}

The fourth result is a large deviation principle for $M_{r,n}(\alpha, \theta)$ in the Ewens-Pitman sampling model \cite{FF15}.

\begin{proposition}\cite[Theorem 1.1]{FF15} \label{ClustersLDPLevel0Frequencies}
For $\alpha \in (0,1)$ and $\alpha > -\theta$, the sequence $\left\{\frac{M_{r,n}(\alpha, \theta)}{n} \right\}$ satisfies a large deviation principle with speed $n$ and rate function $I_r^\alpha(x) = \sup_{\lambda \in \R} \{ \lambda x - \Lambda_{\alpha,r}(\lambda)\}$, 
where
\bal
\Lambda_{\alpha, r}(\lambda) = \begin{cases}
\log\left( 1 + \frac{\alpha \epsilon_0(\lambda)}{1 - r\epsilon_0(\lambda)} \right) & \text{if $\lambda > 0$}, \\
0 & \text{otherwise},
\end{cases}
\eal
with $\epsilon_0(\lambda)$ the unique solution to the equation
\bal
(r - \alpha) \log(1 - (r - \alpha)\epsilon) - r\log(1 - r\epsilon) - \alpha \log(\alpha \epsilon) - \log\left( p_\alpha(r) \frac{1 - e^{-\lambda}}{e^{-\lambda}}  \right) = 0, 
\eal
for $\lambda > 0$ and  $r \geq 1$.
\end{proposition}

The fifth result is a law of the iterated logarithms for $K_n(\alpha, \theta)$ in the Ewens-Pitman sampling model \cite{BF24}.

\begin{proposition}[\cite{BF24}, Theorem 2.4] \label{EPMLIL}
Let $K_n(\alpha, \theta)$ be the number of clusters in a random partition of $[n]$ distributed as $\EPM(\alpha, \theta)$, with $\alpha \in (0,1)$ and $\theta > -\alpha$. Then
\bal
\limsup_{n \to \infty} \left( \frac{K_n - n^\alpha S_{\alpha, \theta}}{\sqrt{2n^\alpha \log\log n}} \right) = \sqrt{S_{\alpha, \theta}}
\eal
almost surely. 
\end{proposition}

The last result is a concentration inequality for $K_n(\alpha, \theta)$ in the Ewens-Pitman sampling model \cite{Per(22)}.

\begin{proposition}\cite[Theorem 3.1]{Per(22)} \label{ClustersConcentration}
Let $K_n(\alpha, \theta)$ be the number of clusters in a random partition of $[n]$ distributed as $\EPM(\alpha, \theta)$, with $\alpha \in (0,1)$ and $\theta > -\alpha$. Then there exists constants $C(\alpha, \theta) > 0$ and $c(\alpha, \theta) > 0$ such that for $\delta < e^{-C(\alpha, \theta)}$, the following holds with probability $\geq 1 - \delta$ for any fixed $n \in \N$:
\bal
\left| \frac{K_n(\alpha, \theta)}{n^\alpha} - S_{\alpha, \theta} \right| \leq \frac{c(\alpha, \theta)[\log \log(n+2) + \log(\delta^{-1})]}{(n + \theta)^{\alpha/2}}.
\eal
\end{proposition}


\section{Almost Sure and $L^p$ Convergence for ${\bf \widetilde M}_{\ell,N}$} \label{Sec:ClustersFrequencyLLN}

We establish the almost sure and $L^p$ convergence for ${\bf \widetilde M}_{\ell,N}(\alpha,\theta_0; \beta,\theta)$, extending  Proposition \ref{HEPLLN}. From the random index representation \eqref{RandomIndex-2}, we start by showing that the random variable $\xi(N)$ in (\ref{totalclusters}), properly scaled, converges in $p$th moment for all integers $p \geq 1$.

\begin{lemma} \label{level1conv}
Let $\{\Pi_{N_i}^{(i)}\}_{1 \leq i \leq d}$ be independent random partitions such that $\Pi_{N_i}^{(i)}$ is a random partition of $[N_i]$ distributed as $\EPM(\beta, \theta)$, for all $1 \leq i \leq d$. 
Let $\xi(N) = \sum_{i=1}^d K_{i,N_i}(\beta, \theta)$, where $\{K_{i,N_i}(\beta, \theta)\}_{1 \leq i \leq d}$ are independent random variables such that $K_{i,N_i}(\beta, \theta)$ is the number of clusters in $\Pi_{N_i}^{(i)}$. 
Suppose $\frac{N_i}{N} = \frac{N_i(N)}{N} \to w_i > 0$ as $N \to \infty$ for all $i \in [d]$.
Then
\begin{enumerate}
\item For $\beta = 0$,
\bal
\frac{\xi(N)}{\log N} &\xrightarrow{a.s.} \theta d \\ 
\E\left[ \left(\frac{\xi(N)}{\log N}\right)^p \right] &\to \E\left[ (\theta d)^p \right]
\eal
as $N \to \infty$ for all integers $p \geq 1$.
\item For $\beta \in (0,1)$,
\bal
\frac{\xi(N)}{N^\beta} &\xrightarrow{a.s.} \eta \\ 
\E\left[ \left(\frac{\xi(N)}{N^\beta}\right)^p \right] &\to \E\left[ \eta^p \right]
\eal
as $N \to \infty$ for all integers $p \geq 1$, where $\eta = \sum_{i=1}^d w_i^\beta S_{\beta,\theta,i}$ such that $\{S_{\beta, \theta, i}\}_{i \in [d]}$ are iid copies of $S_{\beta,\theta}$.
\end{enumerate}
\end{lemma}

\begin{proof}
We only prove Case (2), since Case (1) follows by a similar argument. 
By direct computation and applying (\ref{old-4}),
\bal
\frac{\xi(N)}{N^\beta} = \sum_{i=1}^d \frac{K_{i,N_i}(\beta, \theta)}{N_i^\beta} \left(\frac{N_i}{N}\right)^\beta \xrightarrow{a.s.} \sum_{i=1}^d S_{\beta, \theta, i} w_i^\beta = \eta,
\eal
where $\eta = \sum_{i=1}^d w_i^\beta S_{\beta,\theta,i}$ such that $\{S_{\beta, \theta, i}\}_{i \in [d]}$ are iid copies of $S_{\beta,\theta}$.

Next note that for fixed $p\geq 1$, we have $\left\{ \left( \frac{\xi(N)}{N^\beta}\right)^p \right\} \in L^p$ and $\eta^p \in L^p$, with $\left(\frac{\xi(N)}{N^\beta}\right)^p \xrightarrow{a.s.} \eta^p$. Then
\bal
\E\left( \left( \frac{\xi(N)}{N^\beta}\right)^p \right) &= \E\left[ \left( \sum_{i=1}^d \frac{K_{i, N_i}(\beta, \theta)}{N_i^\beta} \left(\frac{N_i}{N}\right)^\beta \right)^p \right]  \\
&= \sum_{\substack{r_1 + \dotsb + r_d = p \\ r_1, \ldots, r_d \geq 0}} \binom{p}{r_1, \ldots, r_d} \prod_{i=1}^d \left( \frac{N_i}{N} \right)^{\beta r_i} \E\left[\left( \frac{K_{i, N_i}(\beta, \theta)}{N_i^\beta} \right)^{r_i}\right] \\
&\to \sum_{\substack{r_1 + \dotsb + r_d = p \\ r_1, \ldots, r_d \geq 0}} \binom{p}{r_1, \ldots, r_d} \prod_{i=1}^d w_i^{\beta r_i} \E\left[\left( S_{\beta,\theta,i} \right)^{r_i}\right]  \\
&= \E\left[\left( \sum_{i=1}^d w_i^\beta S_{\beta, \theta, i} \right)^p\right] \\
&= \E(\eta^p)
\eal
as $N \to \infty$, where we used the independence of $\{K_{i, N_i}(\beta, \theta)\}_{i \in [d]}$ and the fact that $\E\left[\left( \frac{K_{i, N_i}(\beta, \theta)}{N_i^\beta} \right)^r \right] \to \E\left[\left( S_{\beta,\theta,i} \right)^r \right]$ for all integers $r \geq 1$. 
\end{proof}

\begin{proposition} \label{HEPblockLLN}
For any $1 \leq \ell \leq N$, let ${\bf \widetilde M}_{\ell,N}(\alpha,\theta_0; \beta,\theta)$ be the number of distinct values represented by exactly $\ell$ first-level clusters in a random sample $(X_{i,j})_{1 \leq i \leq d, 1 \leq j \leq N_i}$ of size $N = \sum_{i=1}^d N_i$ from $\HPYP(\alpha, \theta_0, \beta, \theta, \nu, d)$. 
Suppose $\frac{N_i}{N} = \frac{N_i(N)}{N} \to w_i > 0$ as $N \to \infty$ for all $i \in [d]$. Let $p_\alpha(\ell) := \frac{\alpha(1-\alpha)^{(\ell-1)}}{\ell!}$. Then as $N \to \infty$, 
\begin{enumerate}
\item If $\alpha \in (0,1)$ and $\beta = 0$, 
\bal
\frac{{\bf \widetilde M}_{\ell,N}(\alpha, \theta_0; \theta)}{(\log N)^\alpha} \xrightarrow{a.s.} p_\alpha(\ell) S_{\alpha, \theta_0} (\theta d)^\alpha
\eal
\item If $\alpha, \beta \in (0,1)$, 
\bal
\frac{{\bf \widetilde M}_{\ell,N}(\alpha, \theta_0; \beta, \theta)}{N^{\alpha\beta}} \xrightarrow{a.s.} p_\alpha(\ell) S_{\alpha, \theta_0} \eta^\alpha
\eal
\end{enumerate}
where $\eta = \sum_{i=1}^d w_i^\beta S_{\beta,\theta,i}$ such that $\{S_{\beta, \theta, i}\}_{i \in [d]}$ are iid copies of $S_{\beta,\theta}$ independently of $S_{\alpha, \theta_0}$. Moreover the above convergences also hold in $L^p$ for all integers $p \geq 1$.
\end{proposition}

\begin{proof}
We only prove Case (2), as the proof of Case (1) is similar. By virtue of the random index representation (\ref{RandomIndex-2}), 
\bal
\frac{{\bf \widetilde M}_{\ell,N}(\alpha, \theta_0; \beta, \theta)}{N^{\alpha\beta}} \stackrel{a.s.}{=} \frac{M_{\ell, \xi(N)}(\alpha, \theta_0)}{\xi(N)^\alpha} \left(\frac{\xi(N)}{N^\beta}\right)^\alpha.
\eal
Since $\xi(N) \xrightarrow{a.s.} \infty$ as $N \to \infty$, it follows from (\ref{old-5}) that $\frac{M_{\ell, \xi(N)}(\alpha, \theta_0)}{\xi(N)^\alpha} \xrightarrow{a.s.} p_\alpha(\ell) S_{\alpha, \theta_0}$ as $N \to \infty$. Moreover by Lemma \ref{level1conv}, $\frac{\xi(N)}{N^\beta} \xrightarrow{a.s.} \eta = \sum_{i=1}^d w_i^\beta S_{\beta,\theta,i}$, as $N \to \infty$. Thus we have that
\bal
\frac{{\bf \widetilde M}_{\ell,N}(\alpha, \theta_0; \beta, \theta)}{N^{\alpha \beta}} \xrightarrow{a.s.} p_\alpha(\ell) S_{\alpha, \theta_0} \eta^\alpha
\eal
as $N \to \infty$. 

It remains to prove convergence in $L^p$.
Since $\xi(N) \xrightarrow{a.s.} \infty$ as $N \to \infty$, it also follows from (\ref{old-5}) that $\frac{M_{\ell, \xi(N)}(\alpha, \theta_0)}{\xi(N)^\alpha} \stackrel{L^p}{\longrightarrow} p_\alpha(\ell)S_{\alpha, \theta_0}$ as $N \to \infty$, for all integers $p \geq 1$. 
Next we show that $\left(\frac{\xi(N)}{N^\beta}\right)^\alpha \stackrel{L^p}{\longrightarrow} \eta^\alpha$. 
For all $p \geq 1$ and $\alpha \in (0,1)$, we have that 
\bal
\left(\frac{\xi(N)}{N^\beta}\right)^{\alpha p} \leq \left(\frac{\xi(N)}{N^\beta}\right)^p.
\eal
By Lemma \ref{level1conv}, $\E\left( \left( \frac{\xi(N)}{N^\beta}\right)^p \right) \to \E(\eta^p)$ as $N \to \infty$, and by the continuous mapping theorem, $\left(\frac{\xi(N)}{N^\beta}\right)^{\alpha p} \xrightarrow{a.s.} \eta^{\alpha p}$. Therefore applying the generalized dominated convergence theorem yields
\bal
\E\left( \left(\frac{\xi(N)}{N^\beta}\right)^{\alpha p} \right) \to \E(\eta^{\alpha p})
\eal
as $N \to \infty$. Since $\left(\frac{\xi(N)}{N^\beta}\right)^{\alpha} \xrightarrow{a.s.} \eta^{\alpha}$, it follows from the Riesz-Scheffe theorem that $\left(\frac{\xi(N)}{N^\beta}\right)^{\alpha} \stackrel{L^p}{\longrightarrow} \eta^{\alpha}$ for all $p \geq 1$. 

Finally by the Minkowski and Cauchy-Schwarz inequalities, 
\bal
&\E\left( \left| \frac{M_{\ell, \xi(N)}(\alpha, \theta_0)}{\xi(N)^\alpha} \left(\frac{\xi(N)}{N^\beta}\right)^\alpha - p_\alpha(\ell) S_{\alpha, \theta_0} \eta^\alpha \right|^p\right) \\ 
&\qquad = \E\left[ \left| \left( \frac{M_{\ell, \xi(N)}(\alpha, \theta_0)}{\xi(N)^\alpha} \left(\frac{\xi(N)}{N^\beta}\right)^\alpha - \frac{M_{\ell, \xi(N)}(\alpha, \theta_0)}{\xi(N)^\alpha} \eta^\alpha \right) \right. \right. \\
& \qquad \qquad \qquad  \left. \left. + \left( \frac{M_{\ell, \xi(N)}(\alpha, \theta_0)}{\xi(N)^\alpha} \eta^\alpha - p_\alpha(\ell) S_{\alpha, \theta_0} \eta^\alpha \right) \right|^p \right] \\
&\qquad \leq 2^{p-1}\left[ \E\left(\left| \frac{M_{\ell, \xi(N)}(\alpha, \theta_0)}{\xi(N)^\alpha}\left( \left(\frac{\xi(N)}{N^\beta}\right)^\alpha - \eta^\alpha\right) \right|^p \right) \right. \\
&\qquad \qquad \qquad \left. + \E\left(\left| \eta^\alpha \left(\frac{M_{\ell, \xi(N)}(\alpha, \theta_0)}{\xi(N)^\alpha} - p_\alpha(\ell) S_{\alpha, \theta_0} \right) \right|^p \right) \right] \\
&\qquad \leq 2^{p-1}\left[ \sqrt{\E\left[ \left(\frac{M_{\ell, \xi(N)}(\alpha, \theta_0)}{\xi(N)^\alpha}\right)^{2p}\right] \E\left[ \left( \left(\frac{\xi(N)}{N^\beta}\right)^\alpha - \eta^\alpha \right)^{2p}\right]} \right. \\
&\qquad \qquad \qquad \left. + \sqrt{\E\left[ \eta^{2\alpha p}\right] \E\left[ \left( \frac{M_{\ell, \xi(N)}(\alpha, \theta_0)}{\xi(N)^\alpha} - p_\alpha(\ell) S_{\alpha, \theta_0} \right)^{2p}\right]} \right] \\
&\qquad \to 0
\eal
as $N \to \infty$, since $\frac{M_{\ell, \xi(N)}(\alpha, \theta_0)}{\xi(N)^\alpha} \stackrel{L^{2p}}{\longrightarrow} p_\alpha(\ell) S_{\alpha, \theta_0}$ and $\left(\frac{\xi(N)}{N^\beta}\right)^\alpha \stackrel{L^{2p}}{\longrightarrow} \eta^\alpha$ as $N \to \infty$. 
Therefore we conclude that $\frac{{\bf \widetilde M}_{\ell, N}(\alpha, \theta_0; \beta, \theta)}{N^{\alpha \beta}} \stackrel{L^p}{\longrightarrow} p_\alpha(\ell) S_{\alpha, \theta_0} \eta^\alpha $ for all integers $p \geq 1$. 
\end{proof}

\begin{remark}
From Proposition \ref{HEPblockLLN}, as $N \to \infty$
\bal
\frac{{\bf K}_N(\alpha, \theta_0, \beta, \theta)}{N^{\alpha\beta}} = \sum_{\ell=1}^{\xi(N)} \frac{M_{\ell,\xi(N)}(\alpha, \theta_0)}{N^{\alpha\beta}} \xrightarrow{a.s.} S_{\alpha, \theta_0} \eta^\alpha,
\eal
using the fact that $\sum_{\ell=1}^n p_\alpha(\ell) \to 1$, as $n \to \infty$. This convergence also holds in $L^p$. Therefore, Proposition \ref{HEPblockLLN} allows us to recover Proposition \ref{HEPLLN}.
\end{remark}

\begin{remark}\label{estim_alpha}
By Propositions \ref{HEPLLN} and \ref{HEPblockLLN}, 
\bal
\frac{{\bf \widetilde M}_{1,N}(\alpha, \theta_0, \beta, \theta)}{{\bf K}_N(\alpha, \theta_0, \beta, \theta)} \xrightarrow{a.s.} \alpha
\eal
as $N \to \infty$. This convergence implies that $\frac{{\bf \widetilde M}_{1,N}(\alpha, \theta_0, \beta, \theta)}{{\bf K}_N(\alpha, \theta_0, \beta, \theta)} $ is a (strongly) consistent estimator of $\alpha\in(0,1)$. 
\end{remark}

We conclude the present section by establishing some Poisson limit theorems for ${\bf \widetilde M}_{\ell,N}(0, \theta_0, \beta, \theta)$.

\begin{proposition} \label{HEPblockPoissonLimit}
For any $1 \leq \ell \leq N$, let ${\bf \widetilde M}_{\ell,N}(\alpha,\theta_0; \beta,\theta)$ be the number of distinct values represented by exactly $\ell$ first-level clusters in a random sample $(X_{i,j})_{1 \leq i \leq d, 1 \leq j \leq N_i}$ of size $N = \sum_{i=1}^d N_i$ from $\HPYP(\alpha, \theta_0, \beta, \theta, \nu, d)$. 
Then as $N \to \infty$, 
\begin{enumerate}
\item If $\alpha = \beta = 0$, 
\bal
{\bf \widetilde M}_{\ell,N}(\theta_0; \theta) \stackrel{D}{\longrightarrow} Z_\ell
\eal
and 
\bal
({\bf \widetilde M}_{1,N}(\theta_0; \theta), {\bf \widetilde M}_{2,N}(\theta_0; \theta), \ldots) \stackrel{D}{\longrightarrow} (Z_1, Z_2, \ldots)
\eal
where $\{Z_\ell\}$ are independent Poisson random variables such that $Z_\ell \sim \Poi\left( \frac{\theta_0}{\ell} \right)$. 
\item If $\alpha = 0$ and $\beta \in (0,1)$, then the above convergences also hold for ${\bf \widetilde M}_{\ell,N}(\theta_0; \beta, \theta)$. 
\end{enumerate}
\end{proposition}

\begin{proof}
We only prove Case (1), as Case (2) follows by the same argument.  
By the random index representation (\ref{RandomIndex-2}), ${\bf \widetilde M}_{\ell,N}(\theta_0; \theta) \stackrel{a.s.}{=} M_{\ell,\xi(N)}(\theta_0)$. Conditioning on $\xi(N)$ and applying (\ref{old-3}), the single level Poisson limit theorem for $M_{\ell,n}(\theta_0)$, we get that 
\bal
\left(M_{\ell,\xi(N)}(\theta_0) \middle\vert \xi(N)\right) \stackrel{D}{\longrightarrow} Z_\ell \sim \Poi\left( \frac{\theta_0}{\ell} \right)
\eal
as $N \to \infty$. In particular, the conditional characteristic function satisfies
\bal
\E\left[ \exp\left( it M_{\ell,\xi(N)}(\theta_0) \right) \middle\vert \xi(N) \right] \to \E(\exp(itZ_\ell)) = \exp\left( \frac{\theta_0}{\ell} \left( e^{it} - 1 \right) \right)
\eal
as $N \to \infty$. By the dominated convergence theorem, the unconditional characteristic function of $M_{\ell,\xi(N)}(\theta_0)$ satisfies
\bal
&\left| \E(\exp(it M_{\ell,\xi(N)}(\theta_0))) - \exp\left( \frac{\theta_0}{\ell} \left( e^{it} - 1 \right) \right) \right| \\
&= \left| \E\left[ \E(\exp(it M_{\ell,\xi(N)}(\theta_0)) \vert \xi(N)) - \exp\left( \frac{\theta_0}{\ell} \left( e^{it} - 1 \right) \right)\right] \right| \\
&\leq \E\left| \E(\exp(it M_{\ell,\xi(N)}(\theta_0)) \vert \xi(N)) - \exp\left( \frac{\theta_0}{\ell} \left( e^{it} - 1 \right) \right) \right| \\
&\to 0
\eal
as $N \to \infty$, from which it follows that 
\bal
{\bf \widetilde M}_{\ell,N}(\theta_0; \theta) \stackrel{D}{\longrightarrow} Z_\ell.
\eal

The joint convergence result follows by a similar argument of conditioning on $\xi(N)$, then applying the corresponding single level joint convergence result. 
\end{proof}


\section{Gaussian Fluctuations} \label{Sec:GaussianFluctuations}

Motivated by the random index representation \eqref{RandomIndex-2}, we establish a Gaussian fluctuation for the number of clusters in a random partition of random size distributed as $\SSM(q_0)$. 

\begin{theorem} \label{ClustersCLT}
Let $K_n = K_n(q_0)$ be the number of clusters in a random partition, $\Pi_n$, of size $[n]$ distributed as $\SSM(q_0)$. Let $\{\zeta(n)\}$ be a sequence of finite, positive random variables. Suppose that the following conditions hold:
\begin{enumerate}
\item \label{Cond:AsympDiv} {\rm (}Almost sure convergence of $K_n$ and $\zeta(n)${\rm )}. The following almost sure convergences hold as $n \to \infty$,
\bal
\frac{K_n}{\gamma_0(n)} \xrightarrow{a.s.} S \quad \text{and} \quad \frac{\zeta(n)}{\gamma(n)} \xrightarrow{a.s.} \eta
\eal
where $S$ and $\eta$ are finite positive random variables, and $\{ \gamma_0(n) \}_{n \geq 1}, \{\gamma(n)\}_{n \geq 1}$ are two sequences such that $\gamma_0(n), \gamma(n) \to \infty$ as $n \to \infty$ and $\gamma_0(n)$ is regularly varying of the form $\gamma_0(n) = n^\alpha L(n)$, such that $\alpha \in [0,1)$ is the index of regular variation and $L(n)$ is a slowly varying function. 

\item \label{Cond:DistConv} {\rm (}Distributional convergence of $K_n$ and $\zeta(n)${\rm )}. The following distributional convergences hold as $n \to \infty$,
\bal
\sqrt{\gamma_0(n)} \left( \frac{K_n}{\gamma_0(n)} - S \right) &\xrightarrow{D} Y, \\
\sqrt{\gamma(n)} \left( \frac{\zeta(n)}{\gamma(n)} - \eta \right) &\xrightarrow{D} \sqrt{T} N(0,1),
\eal
where $T$ is a finite positive random variable and $Y$ is a random variable. 
\item \label{Cond:ConvRate} {\rm (}Convergence rate assumption{\rm )}. $\sqrt{\gamma_0(\gamma(n))}\left(\frac{L(\zeta(n))}{L(\gamma(n))} - 1 \right) \to 0$ as $n \to \infty$. 
\end{enumerate}
Then 
\bal
\sqrt{\gamma_0(\gamma(n))}\left( \frac{K_{\zeta(n)}}{\gamma_0(\gamma(n))} - S \eta^\alpha \right) \xrightarrow{D} Y\sqrt{\eta^\alpha}  
\eal
as $n \to \infty$. 
\end{theorem}

\begin{proof}
We start with the following decomposition,
\bal
&\sqrt{\gamma_0(\gamma(n))}\left( \frac{K_{\zeta(n)}}{\gamma_0(\gamma(n))} - S \eta^\alpha \right) \\
&\qquad =  \eta^\alpha \sqrt{\gamma_0(\gamma(n))}\left( \frac{K_{\zeta(n)}}{\gamma_0(\zeta(n))} \left[1 + \frac{\gamma_0(\zeta(n))}{\gamma_0(\gamma(n)) \eta^\alpha} - 1 \right] - S \right) \\
&\qquad = \eta^\alpha \sqrt{ \frac{\gamma_0(\gamma(n))}{\gamma_0(\zeta(n))} } \left[ \sqrt{\gamma_0(\zeta(n))}\left( \frac{K_{\zeta(n)}}{\gamma_0(\zeta(n))} - S \right) \right] \\
&\qquad \qquad + \frac{K_{\zeta(n)}}{\gamma_0(\zeta(n))} \sqrt{\gamma_0(\gamma(n))}\left( \frac{\gamma_0(\zeta(n))}{\gamma_0(\gamma(n))} - \eta^\alpha \right)
\eal

Consider the first summand. Since $\zeta(n) \xrightarrow{a.s.} \infty$ as $n \to \infty$ by condition (\ref{Cond:AsympDiv}), conditioning on $\zeta(n)$ and invoking condition (\ref{Cond:DistConv}) gives
\bal
\left[\sqrt{\gamma_0(\zeta(n))}\left( \frac{K_{\zeta(n)}}{\gamma_0(\zeta(n))} - S \right) \middle\vert \zeta(n) \right] \xrightarrow{D} Y. 
\eal
as $n \to \infty$. Thus the conditional characteristic function satisfies
\bal
\E\left[ \exp\left( it \cdot \sqrt{\gamma_0(\zeta(n))}\left( \frac{K_{\zeta(n)}}{\gamma_0(\zeta(n))} - S \right) \right) \middle\vert \zeta(n) \right] \to \E(\exp(itY))
\eal
as $n \to \infty$. By the dominated convergence theorem, the unconditional characteristic function satisfies
\bal
&\left| \E\left[\exp\left( it \cdot \sqrt{\gamma_0(\zeta(n))}\left( \frac{K_{\zeta(n)}}{\gamma_0(\zeta(n))} - S \right) \right) \right] - \E(\exp(itY)) \right| \\
&= \left| \E\left[ \E\left[ \exp\left( it \cdot \sqrt{\gamma_0(\zeta(n))}\left( \frac{K_{\zeta(n)}}{\gamma_0(\zeta(n))} - S \right) \right) \middle\vert \zeta(n) \right] - \E(\exp(itY)) \right] \right| \\
&\leq \E\left| \E\left[ \exp\left( it \cdot \sqrt{\gamma_0(\zeta(n))}\left( \frac{K_{\zeta(n)}}{\gamma_0(\zeta(n))} - S \right) \right) \middle\vert \zeta(n) \right] - \E(\exp(itY)) \right| \\
&\to 0
\eal
as $n \to \infty$. Consequently, the unconditional convergence in distribution also holds,
\bal
\sqrt{\gamma_0(\zeta(n))}\left( \frac{K_{\zeta(n)}}{\gamma_0(\zeta(n))} - S \right)  
\xrightarrow{D} Y. 
\eal
as $n \to \infty$. Now by condition (\ref{Cond:AsympDiv}),
\bal
\eta^\alpha \sqrt{ \frac{\gamma_0(\gamma(n))}{\gamma_0(\zeta(n))} } = \eta^\alpha \sqrt{ \frac{\gamma(n)^\alpha L(\gamma(n))}{\zeta(n)^\alpha L(\zeta(n))} } \xrightarrow{a.s.} \eta^\alpha \sqrt{\frac{1}{\eta^\alpha}} = \sqrt{\eta^\alpha}
\eal
as $n \to \infty$. Thus the first summand converges in distribution as
\bal
\eta^\alpha \sqrt{ \frac{\gamma_0(\gamma(n))}{\gamma_0(\xi(n))} } \left[ \sqrt{\gamma_0(\zeta(n))}\left( \frac{K_{\zeta(n)}}{\gamma_0(\zeta(n))} - S \right) \right] \xrightarrow{D} Y\sqrt{\eta^\alpha}
\eal
as $n \to \infty$. 

Next consider the second summand, which we can further decompose as
\bal
&\frac{K_{\zeta(n)}}{\gamma_0(\zeta(n))} \sqrt{\gamma_0(\gamma(n))}\left( \frac{\gamma_0(\zeta(n))}{\gamma_0(\gamma(n))} - \eta^\alpha \right) \\
&\qquad = \frac{K_{\zeta(n)}}{\gamma_0(\zeta(n))} \sqrt{\gamma_0(\gamma(n))}\left( \left( \frac{\zeta(n)}{\gamma(n)}\right)^\alpha \left( 1 + \frac{L(\zeta(n))}{L(\gamma(n))} - 1\right) - \eta^\alpha \right) \\
&\qquad = \frac{K_{\zeta(n)}}{\gamma_0(\zeta(n))} \left[ \sqrt{\frac{\gamma_0(\gamma(n))}{\gamma(n)}} \sqrt{\gamma(n)} \left( \left( \frac{\zeta(n)}{\gamma(n)}\right)^\alpha - \eta^\alpha \right) \right. \\
&\qquad \qquad \qquad \qquad \qquad \left. + \left( \frac{\zeta(n)}{\gamma(n)}\right)^\alpha \sqrt{\gamma_0(\gamma(n))} \left(\frac{L(\zeta(n))}{L(\gamma(n))} - 1\right) \right].
\eal
By condition (\ref{Cond:AsympDiv}), $\frac{K_{\zeta(n)}}{\gamma_0(\zeta(n))} \xrightarrow{a.s.} S$. Next, note that $\sqrt{\frac{\gamma_0(\gamma(n))}{\gamma(n)}} = \sqrt{\gamma(n)^{\alpha-1} L(\gamma(n))} \to 0$ since $\alpha - 1 < 0$. 
By condition (\ref{Cond:DistConv}) and the delta method, 
\bal
\sqrt{\gamma(n)} \left( \left( \frac{\zeta(n)}{\gamma(n)}\right)^\alpha - \eta^\alpha \right) \xrightarrow{D} N(0,\tau^2)
\eal
for some random variable $\tau$. By condition (\ref{Cond:ConvRate}), 
\bal
\left( \frac{\zeta(n)}{\gamma(n)}\right)^\alpha \sqrt{\gamma_0(\gamma(n))} \left(\frac{L(\zeta(n))}{L(\gamma(n))} - 1\right) \xrightarrow{a.s.} 0.
\eal
Combining the above and applying Slutsky's theorem several times, we get that 
\bal
\frac{K_{\zeta(n)}}{\gamma_0(\zeta(n))} \sqrt{\gamma_0(\gamma(n))}\left( \frac{\gamma_0(\zeta(n))}{\gamma_0(\gamma(n))} - \eta^\alpha \right) \xrightarrow{D} 0. 
\eal
as $n \to \infty$. 

Finally since convergence in distribution to a constant is equivalent to convergence in probability to that same constant, a final application of Slutsky's theorem yields
\bal
\sqrt{\gamma_0(\gamma(n))}\left( \frac{K_{\zeta(n)}}{\gamma_0(\gamma(n))} - S \eta^\alpha \right) \xrightarrow{D} Y \sqrt{\eta^\alpha}  
\eal
as $n \to \infty$. 
\end{proof}

We apply Theorem \ref{ClustersCLT} to establish a Gaussian fluctuation for $K_N(\alpha, \theta_0, \beta, \theta)$. By relying on the random index representation, our strategy consists of applying Theorem~\ref{ClustersCLT} with $\zeta(n)=\xi(n)$.  

\begin{theorem} \label{HEPCLT}
Let ${\bf K}_N(\alpha, \theta_0; \beta, \theta)$ be the number of distinct values in a random sample $(X_{i,j})_{1 \leq i \leq d, 1 \leq j \leq N_i}$ of size $N = \sum_{i=1}^d N_i$ from $\HPYP(\alpha, \theta_0, \beta, \theta, \nu, d)$. 
Let $\xi(N)$ be the random variable defined in (\ref{totalclusters}). Suppose $\frac{N_i}{N} \to w_i > 0$ as $N \to \infty$ for all $i \in [d]$. Moreover, for $\beta \in (0,1)$ suppose that $\sqrt{N^\beta}\left( \left(\frac{N_i}{N w_i}\right)^\beta - 1 \right) \to 0$ as $N \to \infty$ for all $i \in [d]$. Then the following Gaussian fluctuations hold as $N \to \infty$:
\begin{enumerate}
\item $\sqrt{\log \log N} \left( \frac{{\bf K}_N(\theta_0; \theta)}{\log \log N} - \theta_0 \right) \xrightarrow{D} \sqrt{\theta_0} N(0,1)$.
\item $\sqrt{\beta \log N}\left( \frac{{\bf K}_N(\theta_0; \beta, \theta)}{\beta \log N} - \theta_0  \right) \xrightarrow{D} \sqrt{\theta_0} N(0,1)$.
\item $\sqrt{(\log N)^\alpha} \left( \frac{{\bf K}_N(\alpha, \theta_0; \theta)}{(\log N)^\alpha} - S_{\alpha, \theta_0}(\theta d)^\alpha \right) \stackrel{D}{\longrightarrow} \sqrt{S'_{\alpha, \theta_0} (\theta d)^\alpha} N(0,1)$, where $S_{\alpha, \theta_0}'$ is a random variable independent of $N(0,1)$ and equals in distribution to $S_{\alpha, \theta_0}$.
\item $\sqrt{N^{\alpha \beta}}\left( \frac{{\bf K}_N(\alpha, \theta_0, \beta, \theta)}{N^{\alpha \beta}} - S_{\alpha, \theta_0} \eta^\alpha \right) \stackrel{D}{\longrightarrow} \sqrt{S_{\alpha, \theta_0} \eta^\alpha} N(0,1)$, where $\eta = \sum_{i=1}^d w_i^\beta S_{\beta,\theta,i}$ such that $\{S_{\beta, \theta, i}\}_{i \in [d]}$ are iid copies of $S_{\beta,\theta}$, and $S_{\alpha, \theta_0}'$ is a random variable which equals in distribution to $S_{\alpha, \theta_0}$ and is independent of $N(0,1)$ and $\eta$. 
\end{enumerate}
\end{theorem}

\begin{proof}
By virtue of the random index representation (\ref{RandomIndex-1}), it suffices to check the three conditions of Theorem \ref{ClustersCLT} in all four cases. First note that in all cases, $\{\xi(N)\}$ is a sequence of finite positive random variables. 

{\textit{Case (1).}} By the random index representation (\ref{RandomIndex-1}),
\bal
{\bf K}_N(\theta_0; \theta) \stackrel{a.s.}{=} K_{\xi(N)}(\theta_0)
\eal
with $\xi(N) = \sum_{i=1}^d K_{i, N_i}(\theta)$. 
Recall from (\ref{old-1}) that $\frac{K_{n}(\theta_0)}{\log n} \xrightarrow{a.s.} \theta_0$, where the normalizing sequence $\{\log n\}$ is a regularly varying function of the form $n^0 \log n$. Moreover by Lemma \ref{level1conv}, $\frac{\xi(N)}{\log N} \xrightarrow{a.s.} \theta d$. By (\ref{old-2}), 
\bal
\sqrt{\log n} \left( \frac{K_{n}(\theta_0)}{\log n} - \theta_0 \right) \xrightarrow{D} \sqrt{\theta_0} N(0,1). 
\eal
as $n \to \infty$. 

Next we show that $\xi(N)$ also satisfies a Gaussian fluctuation. 
\bal
\sqrt{\log N} \left( \frac{\xi(N)}{\log N} - \theta d \right) &= \sqrt{\log N} \left( \frac{\xi(N)}{\log N} - \sum_{i=1}^d \frac{K_{i, N_i}(\theta)}{\log N_i} + \sum_{i=1}^d \frac{K_{i, N_i}(\theta)}{\log N_i} - \theta d \right) \\
&= \sqrt{\log N} \left( \frac{\xi(N)}{\log N} - \sum_{i=1}^d \frac{K_{i, N_i}(\theta)}{\log N_i} \right) + \sqrt{\log N} \left(\sum_{i=1}^d \frac{K_{i, N_i}(\theta)}{\log N_i} - \theta d \right) \\
&= \sum_{i=1}^d \frac{K_{i, N_i}(\theta)}{\log N_i} \sqrt{\log N} \left( \frac{\log N_i}{\log N} - 1 \right) \\
&\qquad + \sum_{i=1}^d \sqrt{\frac{\log N}{\log N_i}} \sqrt{\log N_i}\left( \frac{K_{i, N_i}(\theta)}{\log N_i} - \theta \right).
\eal
The first term converges to $0$ since $\frac{K_{i, N_i}(\theta)}{\log N_i} \xrightarrow{a.s.} \theta$ and $\sqrt{\log N} \left( \frac{\log N_i}{\log N} - 1 \right) \to 0$. Since 
\bal
\sqrt{\log N_i} \left( \frac{K_{i, N_i}(\theta)}{\log N_i} - \theta \right) \xrightarrow{D} \sqrt{\theta}N(0,1)
\eal
for all $i \in [d]$ and $\{K_{i, N_i}(\theta)\}$ are independent, the following multivariate Gaussian fluctuations hold
\bal
\begin{pmatrix}
\sqrt{\log N_1} \left( \frac{K_{1,N_1}(\theta)}{\log N_1} - \theta \right) \\ \vdots \\ \sqrt{\log N_d} \left( \frac{K_{d,N_d}(\theta)}{\log N_d} - \theta \right)
\end{pmatrix} \xrightarrow{D} 
\begin{pmatrix}
\sqrt{\theta} N(0,1) \\ \vdots \\ \sqrt{\theta} N(0,1)
\end{pmatrix} \sim \sqrt{\theta} N\left( \bm{0}, \bm{I} \right)
\eal
as $N \to \infty$, where $\bm{0}$ is the zero vector and $\bm{I}$ is the $d \times d$ identity matrix. By Slutsky's theorem
\bal
\sum_{i=1}^d \sqrt{\frac{\log N}{\log N_i}} \sqrt{\log N_i}\left( \frac{K_{i,N_i}(\theta)}{\log N_i} - \theta \right) \xrightarrow{D} \sqrt{\theta} d N(0,1). 
\eal
Combining the above and applying Slutsky's theorem,
\bal
\sqrt{\log N} \left( \frac{\xi(N)}{\log N} - \theta d \right) \xrightarrow{D} \sqrt{\theta d} N(0,1). 
\eal

Finally by Lemma \ref{level1conv}, $\sqrt{\log \log N} \left( \frac{\log \xi(N)}{\log\log N} - 1 \right) \to 0$ as $N \to \infty$, so that condition (\ref{Cond:ConvRate}) is satisfied. Therefore all three conditions of Theorem \ref{ClustersCLT} are satisfied, and it follows that 
\bal
\sqrt{\log\log N}\left( \frac{{\bf K}_N(\theta_0; \theta)}{\log\log N} - \theta_0 \right) \stackrel{D}{\longrightarrow} \sqrt{ \theta_0} N(0,1)
\eal
as $N \to \infty$. 

{\textit{Case (2).}} By the random index representation (\ref{RandomIndex-1}),
\bal
{\bf K}_N(\theta_0; \beta, \theta) \stackrel{a.s.}{=} K_{\xi(N)}(\theta_0)  
\eal
with $\xi(N) = \sum_{i=1}^d K_{i, N_i}(\beta, \theta)$. 
Again by (\ref{old-1}), $\frac{K_{n}(\theta_0)}{\log n} \xrightarrow{a.s.} \theta_0$ and by (\ref{old-2}), 
\bal
\sqrt{\log n} \left( \frac{K_{n}(\theta_0)}{\log n} - \theta_0 \right) \xrightarrow{D} \sqrt{\theta_0} N(0,1). 
\eal
Next by Lemma \ref{level1conv}, $\frac{\xi(N)}{N^\beta} \xrightarrow{a.s.} \eta = \sum_{i=1}^d w_i^\beta S_{\beta, \theta, i}$.

Next we show that $\xi(N)$ also satisfies a Gaussian fluctuation. Consider
\bal
\sqrt{N^\beta}\left( \frac{\xi(N)}{N^\beta} - \eta \right) &= \sqrt{N^\beta}\left( \frac{\xi(N)}{N^\beta} - \sum_{i=1}^d \frac{K_{i, N_i}(\beta, \theta) w_i^\beta}{N_i^\beta} + \sum_{i=1}^d \frac{K_{i, N_i}(\beta, \theta) w_i^\beta}{N_i^\beta} - \eta \right) \\
&= \sqrt{N^\beta}\left( \frac{\xi(N)}{N^\beta} - \sum_{i=1}^d \frac{K_{i, N_i}(\beta, \theta) w_i^\beta}{N_i^\beta} \right) + \sqrt{N^\beta}\left(\sum_{i=1}^d \frac{K_{i, N_i}(\beta, \theta) w_i^\beta}{N_i^\beta} - \eta \right) \\
&= \sum_{i=1}^d \frac{K_{i, N_i}(\beta, \theta)}{N_i^\beta} \sqrt{N^\beta}\left( \left(\frac{N_i}{N w_i}\right)^\beta - 1 \right) \\
&\qquad + \sum_{i=1}^d \sqrt{\left( \frac{N}{N_i} \right)^\beta} w_i^\beta \sqrt{N_i^\beta} \left(\frac{K_{i, N_i}(\beta, \theta)}{N_i^\beta} - S_{\beta,\theta, i} \right)
\eal
The first term converges to $0$ as $N \to \infty$ by hypothesis.
For the second term, we apply Theorem \ref{EPCLT} and Slutsky's theorem to get
\bal
\sqrt{\left(\frac{N}{N_i}\right)^\beta} w_i^\beta \sqrt{N_i^\beta} \left( \frac{K_{i,N_i}(\beta, \theta)}{N_i^\beta} - S_{\beta, \theta, i} \right) \stackrel{D}{\longrightarrow} \sqrt{w_i^\beta S_{\beta,\theta,i}'} N(0,1)
\eal
as $N \to \infty$ for all $i \in [d]$, where $\{S_{\beta,\theta,i}'\}_{i \in [d]}$ are iid random variables independent of $N(0,1)$ and equal in distribution to $S_{\beta,\theta}$. Moreover since $\{K_{i,N_i}(\beta, \theta)\}_{i=1}^d$ are independent, the following joint convergence holds
\bal
\begin{pmatrix}
\sqrt{\left(\frac{N}{N_1}\right)^\beta} w_1^\beta \sqrt{N_1^\beta} \left( \frac{K_{1,N_1}(\beta, \theta)}{N_1^\beta} - S_{\beta, \theta, 1} \right) \\
\vdots \\
\sqrt{\left(\frac{N}{N_d}\right)^\beta} w_d^\beta \sqrt{N_d^\beta} \left( \frac{K_{d,N_d}(\beta, \theta)}{N_d^\beta} - S_{\beta, \theta, d} \right) 
\end{pmatrix} \stackrel{D}{\longrightarrow}
\begin{pmatrix}
\sqrt{w_1^\beta S'_{\beta,\theta,1}} N(0,1) \\
\vdots \\
\sqrt{w_d^\beta S'_{\beta,\theta,d}} N(0,1)
\end{pmatrix} \sim N\left( \bm{0}, \bm{\Sigma} \right)
\eal
as $N \to \infty$, where the $d \times d$ covariance matrix $\bm{\Sigma} = (\bm{\Sigma}_{ij})_{1 \leq i,j \leq d}$ is given by 
\bal
\bm{\Sigma}_{ij} = \begin{cases}
w_i^\beta S'_{\beta, \theta, i} & \text{for $i = j$} \\
0 & \text{otherwise}.
\end{cases}
\eal
Applying Slutsky's theorem again gives
\bal
\sum_{i=1}^d \sqrt{\left(\frac{N}{N_i}\right)^\beta} w_i^\beta \sqrt{N_i^\beta} \left( \frac{K_{i,N_i}(\beta,\theta)}{N_i^\beta} - S_{\beta, \theta, i} \right) \stackrel{D}{\longrightarrow} \sum_{i=1}^d \sqrt{w_i^\beta S'_{\beta,\theta, i}} N(0,1)
\eal
as $N \to \infty$. Combining the above and applying Slutsky's theorem once more yields
\bal
\sqrt{N^\beta} \left( \frac{\xi(N)}{N^\beta} - \eta \right) \stackrel{D}{\longrightarrow} \sum_{i=1}^d \sqrt{w_i^\beta S'_{\beta,\theta, i}} N(0,1). 
\eal
Finally $\sqrt{\beta \log N}\left( \frac{\log \xi(N)}{\beta \log N} - 1 \right)$ as $N \to \infty$ by Lemma \ref{level1conv}. 

Therefore all three conditions of Theorem \ref{ClustersCLT} are satisfied and it follows that 
\bal
\sqrt{\beta \log N}\left( \frac{{\bf K}_N(\theta_0; \beta, \theta)}{\beta \log N} - \theta_0  \right) \xrightarrow{D} \sqrt{\theta_0} N(0,1). 
\eal

{\textit{Case (3).}} We have the almost sure representation
\bal
{\bf K}_N(\alpha, \theta_0; \theta) \stackrel{a.s.}{=} K_{\xi(N)}(\alpha, \theta_0) 
\eal
with $\xi(N) = \sum_{i=1}^d K_{i, N_i}(\theta)$. 
By (\ref{old-4}), $\frac{K_{n}(\alpha, \theta_0)}{n^\alpha} \xrightarrow{a.s.} S_{\alpha, \theta_0}$, with regularly varying normalizing sequence $\{n^\alpha\}$. By Proposition~\ref{EPCLT},
\bal
\sqrt{n^\alpha} \left( \frac{K_{n}(\alpha, \theta_0)}{n^\alpha} - S_{\alpha, \theta_0} \right) \xrightarrow{D} \sqrt{S'_{\alpha, \theta_0}} N(0,1),
\eal
where $S_{\alpha, \theta_0}'$ is a random variable independent of $N(0,1)$ and equals in distribution to $S_{\alpha, \theta}$. Moreover by Lemma \ref{level1conv}, $\frac{\xi(N)}{\log N} \xrightarrow{a.s.} \theta d$, and in the proof of Case 1 we showed that $\xi(N)$ satisfies
\bal
\sqrt{\log N} \left( \frac{\xi(N)}{\log N} - \theta d \right) \stackrel{D}{\longrightarrow} \sqrt{\theta d} N(0,1). 
\eal
Finally since our regularly varying sequence is $n^\alpha$, its slowly varying component is $L(n) = 1$ for all $n \geq 1$ and condition (3) is trivially satisfied.

Therefore all three conditions of Theorem \ref{ClustersCLT} are satisfied and it follows that 
\bal
\sqrt{(\log N)^\alpha} \left( \frac{{\bf K}_N(\alpha, \theta_0; \theta)}{(\log N)^\alpha} - S_{\alpha, \theta_0}(\theta d)^\alpha \right) \stackrel{D}{\longrightarrow} \sqrt{S'_{\alpha, \theta_0} (\theta d)^\alpha} N(0,1).
\eal
as $N \to \infty$. 

{\textit{Case (4).}} We have the almost sure representation
\bal
{\bf K}_N(\alpha, \theta_0, \beta, \theta) \stackrel{a.s.}{=} K_{\xi(N)}(\alpha, \theta_0)  
\eal
with $\xi(N) = \sum_{i=1}^d K_{i, N_i}(\beta, \theta)$. Again by (\ref{old-4}), $\frac{K_{n}(\alpha, \theta_0)}{n^\alpha} \xrightarrow{a.s.} S_{\alpha, \theta_0}$ with regularly varying normalizing sequence $\{n^\alpha\}$. By Proposition~\ref{EPCLT},
\bal
\sqrt{n^\alpha} \left( \frac{K_{n}(\alpha, \theta_0)}{n^\alpha} - S_{\alpha, \theta_0} \right) \xrightarrow{D} \sqrt{S'_{\alpha, \theta_0}} N(0,1)
\eal
where $S_{\alpha, \theta_0}'$ is a random variable independent of $N(0,1)$ and equals in distribution to $S_{\alpha, \theta}$.
Moreover by Lemma \ref{level1conv}, $\frac{\xi(N)}{N^\beta} \xrightarrow{a.s.} \eta$, and in the proof of Case 2 we showed 
\bal
\sqrt{N^\beta} \left( \frac{\xi(N)}{N^\beta} - \eta \right) \stackrel{D}{\longrightarrow} \sum_{i=1}^d \sqrt{w_i^\beta S'_{\beta,\theta, i}} N(0,1). 
\eal
Finally as in Case 3, condition (3) is trivially satisfied. 

Therefore all three conditions of Theorem \ref{ClustersCLT} are satisfied and it follows that 
\bal
\sqrt{N^{\alpha \beta}}\left( \frac{{\bf K}_N(\alpha, \theta_0, \beta, \theta)}{N^{\alpha \beta}} - S_{\alpha, \theta_0} \eta^\alpha \right) \stackrel{D}{\longrightarrow} \sqrt{ S'_{\alpha, \theta_0} \eta^\alpha} N(0,1)
\eal
as $N \to \infty$. 
\end{proof}

We conclude the present section with an alternative proof of Proposition~\ref{EPCLT}. Differently from \cite{BF24}, where a martingale construction of $K_{n}(\alpha,\theta)$ is derived to obtain the Gaussian fluctuation of Proposition~\ref{EPCLT}, our proof relies on some direct moment calculations.

\begin{proof}[Proof of Proposition~\ref{EPCLT}]
We use the method of moments to show that
\bal
\E\left[ n^{\frac{r\alpha}{2}}\left( \frac{K_n(\alpha, \theta)}{n^\alpha} - S_{\alpha, \theta} \right)^r \right] \to \E(Z^r)\E(S_{\alpha, \theta}'^{r/2}),
\eal
as $n \to \infty$ for all $r \geq 1$, where $Z \sim N(0,1)$ and is independent of $S_{\alpha, \theta}' \stackrel{D}{=} S_{\alpha, \theta}$. 

We begin by recalling from \cite[Proposition 2]{FLMP09} that the conditional distribution of $S_{\alpha, \theta}$ given $K_n(\alpha, \theta)$ is a {\em beta-tilted $\alpha$-diversity},
\bal
\left( S_{\alpha, \theta} \middle\vert K_n(\alpha, \theta) = k \right) \stackrel{D}{=} B_{k + \frac{\theta}{\alpha}, \frac{n}{\alpha} - k} S_{\alpha, \theta + n},
\eal
where $B_{k + \frac{\theta}{\alpha}, \frac{n}{\alpha} - k} \sim \Beta\left( k + \frac{\theta}{\alpha}, \frac{n}{\alpha} - k \right)$ and is independent of $S_{\alpha, \theta + n}$. 

Conditioning on $K_n(\alpha, \theta)$, we have that
\bal
&\E\left[ n^{\frac{r\alpha}{2}}\left( \frac{K_n(\alpha, \theta)}{n^\alpha} - S_{\alpha, \theta} \right)^r \right] = \E\left( \E\left[ n^{\frac{r\alpha}{2}}\left( \frac{K_n(\alpha, \theta)}{n^\alpha} - S_{\alpha, \theta} \right)^r \middle\vert K_n(\alpha, \theta) \right] \right) \\
&= \E\left( n^\frac{r\alpha}{2} \sum_{j=0}^r \binom{r}{j} (-1)^j \left( \frac{K_n(\alpha, \theta)}{n^\alpha} \right)^{r-j} \left( K_n(\alpha, \theta) + \frac{\theta}{\alpha} \right)^{(j)} \frac{\Gamma(\theta + n)}{\Gamma(\theta + n + j\alpha)} \right) \\
&\sim \E\left( \sum_{j=0}^r (-1)^j \binom{r}{j} \frac{K_n(\alpha, \theta)^{r-j}}{n^{\alpha r/2}}\left( \sum_{\ell = 0}^j {j \brack \ell} K_n(\alpha, \theta)^\ell \right) \right) \\
&= \E\left( \sum_{j=0}^r \frac{K_n(\alpha, \theta)^{r-j}}{n^{\alpha r/2}} \sum_{\ell = j}^r (-1)^\ell \binom{r}{\ell} {\ell \brack \ell - j} \right). 
\eal

First suppose $r$ is even, so that $r = 2m$ for some integer $m \geq 1$. There are three subcases: 
\begin{itemize}
\item If $j \geq m + 1$, then $2m - j \leq m-1$, so that $\E\left( \frac{K_n(\alpha, \theta)^{2m-j}}{n^{\alpha m}} \right) \to 0$ as $n \to \infty$. 
\item If $j \leq m - 1$, then $2m - j \geq m + 1$. Then by applying \cite[Theorem 4]{GKQ15} (making the change of variables $r = 2n+m, j = n, \ell = k$ to go from our notation to theirs), we have that the following identity holds: $\sum_{\ell = j}^{2m} (-1)^\ell \binom{2m}{\ell} {\ell \brack \ell - j} = 0$. 
\item If $j = m$, then $2m - j = m$. Then $\E\left( \frac{K_n(\alpha, \theta)^{m}}{n^{\alpha m}} \right) \to \E(S_{\alpha, \theta}^m)$ as $n \to \infty$ and again by \cite[Theorem 4]{GKQ15} the following identity holds: $\sum_{\ell = m}^{2m} (-1)^\ell \binom{2m}{\ell} {\ell \brack \ell - m} = (2m-1)!!$
\end{itemize}
It follows that for even $r$,
\bal
\E\left[ n^{\frac{r\alpha}{2}}\left( \frac{K_n(\alpha, \theta)}{n^\alpha} - S_{\alpha, \theta} \right)^r \right] \to (r-1)!! \E\left(S_{\alpha, \theta}^{r/2}\right)
\eal
as $n \to \infty$. 

Otherwise suppose $r$ is odd, so that $r = 2m + 1$ for $m \geq 0$. There are two subcases: 
\begin{itemize}
\item If $j \geq m + 1$, then $2m + 1 - j < m$, so that $\E\left( \frac{K_n(\alpha, \theta)^{r-j}}{n^{\alpha r/2}} \right) \to 0$ as $n \to \infty$. 
\item If $j \leq m$, then $2m + 1 - j > m$. Then by  \cite[Theorem 4]{GKQ15}, $\sum_{\ell = j}^{2m+1} (-1)^\ell \binom{2m+1}{\ell} {\ell \brack \ell - j} = 0$. 
\end{itemize}
It follows that for $r$ odd, 
\bal
\E\left[ n^{\frac{r\alpha}{2}}\left( \frac{K_n(\alpha, \theta)}{n^\alpha} - S_{\alpha, \theta} \right)^r \right] \to 0
\eal
as $n \to \infty$. 

Combining the above, we have that
\bal
\E\left[ n^{\frac{r\alpha}{2}}\left( \frac{K_n(\alpha, \theta)}{n^\alpha} - S_{\alpha, \theta} \right)^r \right] \to \E(Z^r) \E(S_{\alpha, \theta}'^{r/2})
\eal
as $n \to \infty$ for all $r \geq 1$, where $Z \sim N(0,1)$ and is independent of $S_{\alpha, \theta}' \stackrel{D}{=} S_{\alpha, \theta}$. Since $Z$ and $S_{\alpha, \theta}$ are independent, and each distribution is determined by its moments, we conclude that 
\bal
\sqrt{n^\alpha}\left( \frac{K_n(\alpha, \theta)}{n^\alpha} - S_{\alpha, \theta} \right) \stackrel{D}{\longrightarrow} \sqrt{S_{\alpha, \theta}'} N(0,1), 
\eal
as $n \to \infty$. 
\end{proof}


\section{Large Deviation Principles} \label{Sec:LDP}

Motivated by the random index representation \eqref{RandomIndex-2}, we establish a large deviation principle for the number of clusters in a random partition of random size distributed as $\SSM(q_0)$.

\begin{theorem} \label{ClustersLDP}
Let $K_n = K_n(q_0)$ be the number of clusters in a random partition of size $n$ distributed as $\SSM(q_0)$. Let $\{\zeta(n)\}$ be a sequence of finite, positive random variables. Suppose that the following conditions hold:
\begin{enumerate}
\item \label{Cond:Level0LDP} Let $\{a_0(n)\}$ be a sequence of positive numbers such that $a_0(n) \to \infty$ as $n \to \infty$. Suppose that
\bal
\Lambda_0(\lambda) = \lim_{n \to \infty} \frac{1}{a_0(n)} \log \E\left[ \exp(\lambda K_n) \right]
\eal
exists for all $\lambda \in \R$. 
\item \label{Cond:LimitExistence} Suppose there exists a sequence $\{a(n)\}$ such that the limit 
\bal
\Lambda_1(\lambda) = \frac{1}{a(n)} \log\E\left[ \exp(\lambda a_0(\zeta(n))) \right]
\eal
exists for all $\lambda \in \R$, and is continuous.
\item The function $\Lambda_1(\Lambda_0(\lambda))$ is essentially smooth.
\end{enumerate}
Then the family $\left\{ \frac{K_{\zeta(n)}}{a(n)} \right\}$ satisfies a large deviation principle with speed $a(n)$ and rate function $I(x) =  \sup_{\lambda \in \R} \{ \lambda x - \Lambda(\lambda) \}$, where $\Lambda(\lambda) := \Lambda_1(\Lambda_0(\lambda))$.
\end{theorem}

\begin{proof}
Fix $\lambda \in \R$. Since the limit 
\bal
\Lambda_0(\lambda) = \lim_{n \to \infty} \frac{1}{a_0(n)} \log \E\left[ \exp(\lambda K_n) \right]
\eal
exists for all $\lambda \in \R$ by Condition \ref{Cond:Level0LDP}, then for $\epsilon > 0$ there exists $M$ sufficiently large such that
\bal
\left| \frac{1}{a_0(n)} \log \E\left[ \exp(\lambda K_n) \right] - \Lambda_0(\lambda) \right| < \epsilon
\eal
for all $n \geq M$. 

Now consider the cumulant generating function of $K_{\zeta(n)}$. Conditioning on $\zeta(n)$, we get
\bal
&\log \E\left[ \exp(\lambda K_{\zeta(n)}) \right] = \log \E\left( \E\left[ \exp(\lambda K_{\zeta(n)}) \mid \zeta(n)\right] \right) \\
&= \log \left( \E\left( \E\left[ \exp(\lambda K_{\zeta(n)}) \mid \xi(n)\right] \I_{\{\zeta(n) \leq M\}} \right) + \E\left( \E\left[ \exp(\lambda K_{\zeta(n)}) \mid \zeta(n)\right] \I_{\{\zeta(n) > M\}} \right) \right).
\eal
Observe that the first expectation term inside the logarithm satisfies
\bal
e^{-|\lambda M|} P(\zeta(n) \leq M) \leq \E\left( \E\left[ \exp(\lambda K_{\zeta(n)}) \mid \zeta(n)\right] \I_{\{\zeta(n) \leq M\}} \right) \leq e^{|\lambda M|} P(\zeta(n) \leq M).
\eal
Since $\zeta(n) \xrightarrow{a.s.} \infty$, it follows that $P(\zeta(n) \leq M) \to 0$ as $n \to \infty$. Consequently,
\bal
\E\left( \E\left[ \exp(\lambda K_{\zeta(n)}) \mid \zeta(n)\right] \I_{\{\zeta(n) \leq M\}} \right) \to 0
\eal
as $n \to \infty$.  Noting that 
\bal
 &\E\left( \E\left[ \exp(\lambda K_{\zeta(n)}) \mid \zeta(n)\right] \I_{\{\zeta(n) > M\}}\right)\\
 &\leq \E\left( \E\left[ \exp(\lambda K_{\zeta(n)}) \mid \zeta(n)\right]\right)= \E\left[ \exp(\lambda K_{\zeta(n)}\right]\\
 & \leq 2\max\bigg\{\E\left( \E\left[ \exp(\lambda K_{\zeta(n)}) \mid \zeta(n)\right] \I_{\{\zeta(n) > M\}}\right), \E\left( \E\left[ \exp(\lambda K_{\zeta(n)}) \mid \zeta(n)\right] \I_{\{\zeta(n) \leq M\}}\right)\bigg\},
 \eal
it follows that 
\bal
\lim_{n \to \infty}\frac{1}{a(n)}\log\E\left[ \exp(\lambda K_{\zeta(n)})\right]=\lim_{n \to \infty}\frac{1}{a(n)}\log\E\left( \E\left[ \exp(\lambda K_{\zeta(n)}) \mid \zeta(n)\right] \I_{\{\zeta(n) > M\}}\right).
\eal

On the other hand, 
\bal
&\E\left( \exp\left[ (\Lambda_0(\lambda) - \epsilon) a_0(\zeta(n)) \right] \I_{\{\zeta(n) > M\}} \right) \\
&\qquad \leq \E\left( \E\left[ \exp(\lambda K_{\zeta(n)}) \mid \zeta(n)\right] \I_{\{\zeta(n) > M\}} \right) \\
&\qquad  = \E\left( \exp\left( a_0(\zeta(n)) \frac{1}{a_0(\zeta(n))} \log\E\left[ \exp(\lambda K_{\zeta(n)}) \mid \zeta(n)\right] \right) \I_{\{\zeta(n) > M\}}\right) \\
&\qquad\qquad  \leq \E\left( \exp\left[ (\Lambda_0(\lambda) + \epsilon) a_0(\zeta(n)) \right] \I_{\{\zeta(n) > M\}} \right). 
\eal
Observe that both $\E\left( \exp\left[ (\Lambda_0(\lambda) + \epsilon) a_0(\zeta(n)) \right] \I_{\{\zeta(n) \leq M\}} \right)$ and $\E\left( \E\left[ \exp((\Lambda_0(\lambda) -\epsilon)a_0{\zeta(n)}) \mid \zeta(n)\right] \I_{\{\zeta(n) \leq M\}}\right)$ converge to zero as $n \to \infty$, it follow that
\bal
\lim_{n \to \infty} \frac{1}{a(n)}\log\E\left( \exp\left[ (\Lambda_0(\lambda) + \epsilon) a_0(\zeta(n)) \right] \I_{\{\zeta(n) > M\}} \right) &= \lim_{n \to \infty} \frac{1}{a(n)}\log\E\left( \exp\left[ (\Lambda_0(\lambda) + \epsilon) a_0(\zeta(n)) \right] \right)
\eal
and 
\bal
\lim_{n \to \infty}\frac{1}{a(n)}\log\E\left( \exp\left[ (\Lambda_0(\lambda) - \epsilon) a_0(\zeta(n)) \right] \I_{\{\zeta(n) > M\}} \right) = \lim_{n \to \infty} \frac{1}{a(n)}\log\E\left( \exp\left[ (\Lambda_0(\lambda) - \epsilon) a_0(\zeta(n)) \right] \right). 
\eal

Putting things together, we obtain that
\bal
&\lim_{n \to \infty} \frac{1}{a(n)}\log\E\left( \exp\left[ (\Lambda_0(\lambda) - \epsilon) a_0(\zeta(n)) \right] \right) \\
&\qquad \leq \lim_{n \to \infty}\frac{1}{a(n)} \log\E\left[ \exp(\lambda K_{\zeta(n)})\right]  \\
&\qquad \qquad \leq \lim_{n \to \infty} \frac{1}{a(n)}\log\E\left( \exp\left[ (\Lambda_0(\lambda) + \epsilon) a_0(\zeta(n)) \right] \right).
\eal
Combining the above and sending $\epsilon \to 0$, it follows from Condition \ref{Cond:LimitExistence} that the limit 
\bal
\Lambda(\lambda) &:= \lim_{n \to \infty} \frac{1}{a(n)} \log \E\left[ \exp(\lambda K_{\zeta(n)}) \right] \\
&= \lim_{n \to \infty} \frac{1}{a(n)} \log \E\left[\exp\left(\Lambda_0(\lambda) a_0(\zeta(n))\right) \right] \\
&= \Lambda_1(\Lambda_0(\lambda))
\eal
exists for all $\lambda \in \R$. Since $\Lambda(\lambda)$ is essentially smooth, it follows from the G\"{a}rtner-Ellis theorem that $\left\{ \frac{K_{\zeta(n)}}{a(n)} \right\}$ satisfies a large deviation principle with speed $a(n)$ and rate function $I(x) =  \sup_{\lambda \in \R} \{ \lambda x - \Lambda(\lambda) \}$ with $\Lambda(\lambda)$ defined as above.  
\end{proof}

We apply Theorem~\ref{ClustersLDP} to establish a large deviation principles for ${\bf K}_N(\alpha, \theta_0; \beta, \theta)$ and ${\bf M}_{r,N}(\alpha, \theta_0; \beta, \theta)$. As for Gaussian fluctuations, our strategy applies Theorem~\ref{ClustersLDP} for $\zeta(n)=\xi(n)$. 

\begin{lemma} \label{LDPsumsinglelevel}
Let $\xi(N) = \sum_{i=1}^d K_{i, N_i}(\beta, \theta)$ be the random variable defined in (\ref{totalclusters}).
\begin{enumerate}
\item For $\beta = 0$, the family $\left\{ \frac{\xi(N)}{\log N} \right\}$ satisfies a large deviation principle with speed $\log N$ and rate function
 \bal
I(x) = \begin{cases}
x\log\left( \frac{x}{\theta d}\right) - x + \theta d, & x \geq 0, \\
\infty, & x < 0
\end{cases}
\eal
\item For $\beta \in (0,1)$, the family $\left\{ \frac{\xi(N)}{N} \right\}$ satisfies a large deviation principle with speed $N$ and rate function
\bal
I(x) = \sup_{\lambda \in \R}\left\{ \lambda x - \Lambda_\beta(\lambda) \sum_{i=1}^d w_i\right\}.
\eal
\end{enumerate}
\end{lemma}

\begin{proof}
The proof follows from Proposition~\ref{ClustersLDPLevel0} and the contraction principle. 
\end{proof}

\begin{theorem}\label{HEPLDP}
Let ${\bf K}_N(\alpha, \theta_0; \beta, \theta)$ be the number of distinct values in a random sample $(X_{i,j})_{1 \leq i \leq d, 1 \leq j \leq N_i}$ of size $N = \sum_{i=1}^d N_i$ from $\HPYP(\alpha, \theta_0, \beta, \theta, \nu, d)$. 
Suppose $\frac{N_i}{N} = \frac{N_i(N)}{N}\to w_i > 0$ as $N \to \infty$ for all $i \in [d]$. The following hold:
\begin{enumerate}
\item the family $\left\{\frac{{\bf K}_N(\theta_0; \theta)}{\log\log N} \right\}_{N \geq 1}$ satisfies a large deviation principle with speed $\log\log N$ and rate function 
\bal
I_1(x) = \begin{cases}
x\log\left( \frac{x}{\theta_0}\right) - x + \theta_0, & x \geq 0, \\
\infty, & x < 0.
\end{cases}
\eal
\item the family $\left\{\frac{{\bf K}_N(\theta_0; \beta, \theta)}{\log N} \right\}_{N \geq 1}$ satisfies a large deviation principle with speed $\log N$ and rate function  
\bal
I_2(x) = \beta I_1\left( \frac{x}{\beta} \right).
\eal
\item the family $\left\{\frac{{\bf K}_N(\alpha, \theta_0; \theta)}{\log N} \right\}_{N \geq 1}$ satisfies a large deviation principle with speed $\log N$ and rate function  
\bal
I_3(x) = \sup_{\lambda \in \R}\left\{ \lambda x - \theta d \left(e^{\Lambda_\alpha(\lambda)} - 1 \right) \right\}, 
\eal
where $\Lambda_{\alpha}(\lambda)$ is given in Proposition~\ref{ClustersLDPLevel0}.
\item the family $\left\{\frac{{\bf K}_N(\alpha, \theta_0, \beta, \theta)}{N} \right\}_{N \geq 1}$ satisfies a large deviation principle with speed $N$ and rate function
\bal
I_4(x) = \sup_{\lambda \in \R}\left\{ \lambda x - \Lambda_\beta(\Lambda_\alpha(\lambda)) \sum_{i=1}^d w_i \right\}. 
\eal
\end{enumerate}
\end{theorem}

\begin{proof}
By virtue of the random index representation, it suffices to check the conditions of Theorem \ref{ClustersLDP} and identify the speed, $a(n)$, in the four different cases. 

{\textit{Case 1.}} By Proposition \ref{ClustersLDPLevel0}, $\left\{ \frac{K_{0,n}(\theta_0)}{\log n} \right\}$ satisfies a large deviation principle with speed $\log n$ and rate function $I_\theta(x) = \sup_{\lambda \in \R} \{ \lambda x - \Lambda_0(\lambda)\}$, where $\Lambda_0(\lambda) = \theta_0(e^\lambda - 1)$. 

Next for all $\lambda \in \R$, consider 
\bal
\frac{1}{\log\log N} \log\E\left[ \exp\left( \lambda \log \xi(N) \right) \right] &=  \frac{1}{\log\log N} \log\E\left[ \xi(N)^\lambda \right] \\ &= \frac{1}{\log\log N} \log\E\left[ (\log N)^\lambda \left(\frac{\xi(N)}{\log N} \right)^\lambda \right] \\ 
&=  \frac{\lambda \log\log N + \log\E\left[ \left(\frac{\xi(N)}{\log N} \right)^\lambda \right]}{\log\log N} \\ &= \lambda + \frac{\log\E\left[ \left(\frac{\xi(N)}{\log N} \right)^\lambda \right]}{\log\log N}. 
\eal
If $\lambda \leq 0$, then by the dominated convergence theorem and the fact that $\frac{\xi(N)}{\log N} \xrightarrow{a.s.} \theta d$ from Lemma \ref{level1conv}, we get
\bal
\frac{\log\E\left[ \left(\frac{\xi(N)}{\log N} \right)^\lambda \right]}{\log\log N} \to 0
\eal
as $N \to \infty$. Otherwise suppose $\lambda > 0$. Then 
\bal
\frac{\log\E\left[ \left(\frac{\xi(N)}{\log N} \right)^\lambda \right]}{\log\log N} \leq \frac{\log\E\left[ \left(\frac{\xi(N)}{\log N} \right)^{\lceil \lambda \rceil} \right]}{\log\log N} \sim \frac{\log\E\left[ \left( \theta d \right)^{\lceil \lambda \rceil} \right]}{\log\log N}  \to 0
\eal
as $N \to \infty$ by Lemma \ref{level1conv}. 
In any case, the limit 
\bal
\Lambda_1(\lambda) = \lim_{N \to \infty} \frac{1}{\log\log N} \log\E\left[ \exp\left( \lambda \log \xi(N) \right) \right] = \lambda
\eal
exists for all $\lambda \in \R$ which implies that  $\Lambda_1(\Lambda_0(\lambda)) = \Lambda_0(\lambda)=\theta_0(e^{\lambda}-1)$.  
Therefore all conditions of Theorem \ref{ClustersLDP} are satisfied. Thus $\left\{ \frac{{\bf K}_N(\theta_0; \theta)}{\log\log N} \right\}$ satisfies a large deviation principle with speed $\log\log N$ and rate function 
\bal
I_1(x) = \sup_{\lambda \in \R}\{\lambda x - \Lambda_0(\lambda)\} = \begin{cases}
x\log\left( \frac{x}{\theta_0}\right) - x + \theta_0, & x \geq 0, \\
\infty, & x < 0.
\end{cases}
\eal

{\textit{Case 2.}} This is similar to Case 1, except in this case we have that for all $\lambda \in \R$, 
\bal
\frac{1}{\log N} \log\E\left[ \exp\left( \lambda \log \xi(N) \right) \right] &=  \frac{1}{\log N} \log\E\left[ \xi(N)^\lambda \right]  = \frac{1}{\log N} \log\E\left[ N^{\beta\lambda} \left(\frac{\xi(N)}{N^\beta} \right)^\lambda \right] \\ 
&=  \frac{\beta \lambda \log N + \log\E\left[ \left(\frac{\xi(N)}{N^\beta} \right)^\lambda \right]}{\log N} = \beta\lambda + \frac{\log\E\left[ \left(\frac{\xi(N)}{N^\beta} \right)^\lambda \right]}{\log N}. 
\eal
If $\lambda \leq 0$, then by the dominated convergence theorem and the fact that $\frac{\xi(N)}{N^\beta} \xrightarrow{a.s.} \sum_{i=1}^d w_i^\beta S_{\beta, \theta, i}$ from Lemma \ref{level1conv}, we have that 
\bal
\frac{\log\E\left[ \left(\frac{\xi(N)}{N^\beta} \right)^\lambda \right]}{\log N} \to 0
\eal
as $N \to \infty$. Otherwise suppose $\lambda > 0$. Then 
\bal
\frac{\log\E\left[ \left(\frac{\xi(N)}{N^\beta} \right)^\lambda \right]}{\log N} \leq \frac{\log\E\left[ \left(\frac{\xi(N)}{N^\beta} \right)^{\lceil \lambda \rceil} \right]}{\log N} \sim \frac{\log\E\left[ \left( \sum_{i=1}^d w_i^\beta S_{\beta, \theta, i} \right)^{\lceil \lambda \rceil} \right]}{\log N} \to 0
\eal
as $N \to \infty$ by Lemma \ref{level1conv}. 
In any case, the limit 
\bal
\Lambda_1(\lambda) = \lim_{N \to \infty} \frac{1}{\log N} \log\E\left[ \exp\left( \lambda \log \xi(N) \right) \right] = \beta\lambda
\eal
exists for all $\lambda \in \R$. 

Therefore all conditions of Theorem \ref{ClustersLDP} are satisfied, and  $\left\{ \frac{{\bf K}_N(\theta_0; \beta, \theta)}{\log N} \right\}$ satisfies a large deviation principle with speed $\log N$ and rate function 
\bal
I_2(x) = \sup_{\lambda \in \R}\{\lambda x - \beta \Lambda_0(\lambda)\} = \beta I_1\left( \frac{x}{\beta} \right).
\eal

{\textit{Case 3.}} By Proposition \ref{ClustersLDPLevel0}, $\left\{ \frac{K_{n}(\alpha, \theta_0)}{n} \right\}$ satisfies a large deviation principle with speed $n$ and rate function $I_\alpha(x) = \sup_{\lambda \in \R} \{ \lambda x - \Lambda_\alpha(\lambda)\}$. 
On the other hand by Lemma \ref{LDPsumsinglelevel}, the limit
\bal
\Lambda_1(\lambda) = \lim_{N \to \infty} \frac{1}{\log N} \log\E\left[ \exp\left( \lambda \xi(N) \right) \right] = \theta d (e^\lambda - 1)
\eal
exists for all $\lambda \in \R$. Since  $\Lambda_1(\Lambda_\alpha(\lambda)) = \theta d (e^{\Lambda_\alpha(\lambda)} - 1)$ is essentially smooth, all  conditions of Theorem \ref{ClustersLDP} are satisfied.  Thus $\left\{ \frac{{\bf K}_N(\alpha, \theta_0; \theta)}{\log N} \right\}$ satisfies a large deviation principle with speed $\log N$ and rate function  
\bal
I_3(x) = \sup_{\lambda \in \R} \left\{ \lambda x - \theta d \left(e^{\Lambda_\alpha(\lambda)} - 1 \right) \right\}. 
\eal

{\textit{Case 4.}} This is the same as Case 3, except that the limit is
\bal
\Lambda_1(\lambda) = \lim_{N \to \infty} \frac{1}{N} \log\E\left[ \exp\left( \lambda \xi(N) \right) \right] = \Lambda_\beta(\lambda) \sum_{i=1}^d w_i,
\eal
which exists for all $\lambda \in \R$ by Lemma \ref{LDPsumsinglelevel}. 

As before  all conditions of Theorem \ref{ClustersLDP} are satisfied with  $$\Lambda_1(\Lambda_\alpha(\lambda)) = \Lambda_\beta(\Lambda_\alpha(\lambda)) \sum_{i=1}^d w_i.$$ Thus $\left\{ \frac{{\bf K}_N(\alpha, \theta_0; \beta, \theta)}{N} \right\}$ satisfies a large deviation principle with speed $N$ and rate function 
\bal
I_4(x) &= \sup_{\lambda \in \R}\left\{\lambda x - \Lambda_\beta(\Lambda_\alpha(\lambda)) \sum_{i=1}^d w_i \right\}. \qedhere
\eal
\end{proof}

We conclude the present section with a large deviation principle for ${\bf \widetilde M}_{\ell,N}(\alpha, \theta_0; \beta, \theta)$, for $1 \leq \ell \leq N$.

\begin{theorem}\label{HEPblockLDP}
For any $1 \leq \ell \leq N$, let ${\bf \widetilde M}_{\ell,N}(\alpha, \theta_0; \beta, \theta)$ be the number of distinct values represented by exactly $\ell$ first-level clusters in a random sample $(X_{i,j})_{1 \leq i \leq d, 1 \leq j \leq N_i}$ of size $N = \sum_{i=1}^d N_i$ from $\HPYP(\alpha, \theta_0, \beta, \theta, \nu, d)$. 
Suppose $\frac{N_i}{N} = \frac{N_i(N)}{N}\to w_i > 0$ as $N \to \infty$ for all $i \in [d]$. The following hold:
\begin{enumerate}
\item the family $\left\{\frac{{\bf \widetilde M}_{\ell,N}(\alpha, \theta_0; \theta)}{\log N} \right\}_{N \geq 1}$ satisfies a large deviation principle with speed $\log N$ and rate function  
\bal
I(x) = \sup_{\lambda \in \R}\left\{ \lambda x - \theta d \left(e^{\Lambda_{\alpha, \ell}(\lambda)} - 1 \right) \right\}. 
\eal
\item the family $\left\{\frac{{\bf \widetilde M}_{\ell,N}(\alpha, \theta_0; \beta, \theta)}{N} \right\}_{N \geq 1}$ satisfies a large deviation principle with speed $N$ and rate function
\bal
I(x) = \sup_{\lambda \in \R}\left\{ \lambda x - \Lambda_\beta(\Lambda_{\alpha, \ell}(\lambda)) \sum_{i=1}^d w_i \right\}. 
\eal
\end{enumerate}
\end{theorem}

\begin{proof}
The proof uses Theorem \ref{ClustersLDP} and Proposition \ref{ClustersLDPLevel0Frequencies}, and is similar to Cases 3 and 4 in the proof of Theorem \ref{HEPLDP}.
\end{proof}


\section{Law of the Iterated Logarithm} \label{Sec:LIL}

In this section, we establish the law of the iterated logarithm for ${\bf K}_N(\alpha, \theta_0, \beta, \theta)$.

\begin{theorem} \label{HEPMLIL}
Let ${\bf K}_N(\alpha, \theta_0; \beta, \theta)$ be the number of distinct values in a random sample $(X_{i,j})_{1 \leq i \leq d, 1 \leq j \leq N_i}$ of size $N = \sum_{i=1}^d N_i$ from $\HPYP(\alpha, \theta_0, \beta, \theta, \nu, d)$. 
Suppose $\frac{N_i}{N} = \frac{N_i(N)}{N}\to w_i > 0$ as $N \to \infty$ for all $i \in [d]$.
\begin{enumerate}
\item Let $\alpha = 0$, $\theta_0 > 0$, $\beta = 0$, and $\theta > 0$. Then
\bal
\limsup_{N \to \infty} \left( \frac{{\bf K}_N(\theta_0; \theta) - \theta_0 \log\log N}{\sqrt{2\theta_0\log\log N \log\log(\log\log N)}} \right) = 1 \quad \text{a.s.}
\eal
\item Let $\alpha = 0$, $\theta_0 > 0$, $\beta \in (0,1)$, and $\theta > -\beta$. Then
\bal
\limsup_{N \to \infty} \left( \frac{{\bf K}_N(\theta_0; \beta, \theta) - \theta_0 \beta \log N}{\sqrt{2\theta_0 \beta \log N \log\log (\log N)}} \right) = 1 \quad \text{a.s.}
\eal
\item Let $\alpha \in (0, 1)$, $\theta_0 > -\alpha$, $\beta = 0$, and $\theta > 0$. Then
\bal
\limsup_{N \to \infty}\left( \frac{{\bf K}_N(\alpha, \theta_0; \theta) - (\log N)^\alpha S_{\alpha, \theta_0} (\theta d)^\alpha}{\sqrt{2(\log N)^\alpha \log\log (\log N)}} \right) = \sqrt{S_{\alpha, \theta_0} (\theta d)^\alpha} \quad \text{a.s.}
\eal
\item Let $\alpha \in (0, 1)$, $\theta_0 > -\alpha$, $\beta \in (0,1)$, and $\theta > -\beta$. Then
\bal
\limsup_{N \to \infty} \left( \frac{{\bf K}_N(\alpha, \theta_0; \beta, \theta) - N^{\alpha\beta} S_{\alpha, \theta_0} \eta^\alpha}{\sqrt{2N^{\alpha\beta} \log\log N}} \right) = \sqrt{S_{\alpha, \theta_0} \eta^\alpha} \quad \text{a.s.}
\eal
where $\eta = \sum_{i=1}^d w_i^\beta S_{\beta,\theta,i}$ such that $\{S_{\beta, \theta, i}\}_{i \in [d]}$ are iid copies of $S_{\beta,\theta}$. 
\end{enumerate}
\end{theorem}

\begin{proof}
We only prove Case (1) (respectively, Case (4)), since Case (2) (respectively, Case (3)) follows by a similar argument. 

We begin with the proof of Case (1). Recall that we can write $K_n(\theta_0)$ as a sum of independent Bernoulli random variables \cite[Chapter 4]{Pit06} such that $K_n(\theta_0) = \sum_{k = 1}^n X_k$, where $X_k \sim \Ber\left( \frac{\theta_0}{\theta_0 + k -1} \right)$. By the law of the iterated logarithm for sums of independent random variables (see for example \cite{Kol29}), we can obtain the following law of the iterated logarithm for $K_n(\theta_0)$,
\bal
\limsup_{n \to \infty} \frac{K_n(\theta_0) - \theta_0 \log n}{\sqrt{2\theta_0 \log n \log\log (\log n)}} = 1 \quad \text{a.s.}
\eal
By the random index representation (\ref{RandomIndex-1}), consider the following decomposition
\bal
&\frac{{\bf K}_N(\theta_0; \theta) - \theta_0 \log\log N}{\sqrt{2\theta_0\log\log N \log\log(\log\log N)}} \\
&\quad = \frac{K_{\xi(N)}(\theta_0) - \theta_0 \log \xi(N)}{\sqrt{2\theta_0 \log \xi(N) \log\log (\log \xi(N))}} \sqrt{\frac{2\theta_0 \log \xi(N) \log\log (\log \xi(N))}{2\theta_0\log\log N \log\log(\log\log N)}} \\
&\qquad + \frac{\theta_0\log \xi(N) - \theta_0 \log\log N}{\sqrt{2\theta_0\log\log N \log\log(\log\log N)}}.
\eal
Note that the second term converges to $1$ almost surely as $N \to \infty$. On the other hand, using the fact that $\xi(N) \to \infty$ a.s. as $N \to \infty$, Lemma \ref{level1conv}, the law of the iterated logarithm for $K_n(\theta_0)$, and the fact that $\log n$ is slowly varying, we get
\bal
\limsup_{N \to \infty} \frac{K_{\xi(N)}(\theta_0) - \theta_0 \log \xi(N)}{\sqrt{2\theta_0 \log \xi(N) \log\log (\log \xi(N))}} \sqrt{\frac{2\theta_0 \log \xi(N) \log\log (\log \xi(N))}{2\theta_0\log\log N \log\log(\log\log N)}} = 1 \quad \text{a.s.}
\eal
Combining the above proves Case (1). 

Next we prove Case (4). By virtue of the random index representation (\ref{RandomIndex-1}), consider the following decomposition
\bal
&\frac{{\bf K}_N(\alpha, \theta_0; \beta, \theta) - N^{\alpha\beta} S_{\alpha, \theta_0} \eta^\alpha}{\sqrt{2N^{\alpha\beta} \log\log N}} \\
&\quad = \frac{K_{\xi(N)}(\alpha, \theta_0) - \xi(N)^\alpha S_{\alpha, \theta_0}}{\sqrt{2\xi(N)^\alpha \log\log(\xi(N))}} \sqrt{\frac{2\xi(N)^\alpha \log\log(\xi(N))}{2N^{\alpha\beta} \log\log N}}  \\
&\qquad + \frac{S_{\alpha, \theta_0}}{\sqrt{2}}\frac{\sqrt{N^{\alpha\beta}} \left( \left(\frac{\xi(N)}{N^{\beta}}\right)^\alpha - \eta^\alpha \right)}{\sqrt{\log\log N^\beta}}. 
\eal
By Lemma \ref{level1conv} and the fact that $\log$ is a slowly varying function, 
\bal
\sqrt{\frac{2\xi(N)^\alpha \log\log(\xi(N))}{2N^{\alpha\beta} \log\log N}} = \sqrt{\frac{2\xi(N)^\alpha \log\log(\xi(N)) \log\log(N^\beta)}{2N^{\alpha\beta} \log\log (N^\beta) \log\log N}} \xrightarrow{a.s.} \sqrt{\eta^\alpha}
\eal
as $N \to \infty$. Moreover since $\xi(N) \to \infty$ almost surely as $N \to \infty$, by Proposition \ref{EPMLIL},
\bal
\limsup_{N \to \infty} \frac{K_{\xi(N)}(\alpha, \theta_0) - \xi(N)^\alpha S_{\alpha, \theta_0}}{\sqrt{2\xi(N)^\alpha \log\log(\xi(N))}} = \sqrt{S_{\alpha, \theta_0}}.
\eal

It remains to show that the second summand converges to $0$ almost surely as $N \to \infty$. By taking the first order Taylor expansion of $\left(\frac{\xi(N)}{N^{\beta}}\right)^\alpha$, we can write
\bal
\frac{\sqrt{N^{\alpha\beta}} \left( \left(\frac{\xi(N)}{N^{\beta}}\right)^\alpha - \eta^\alpha \right)}{\sqrt{\log\log N}} \stackrel{a.s.}{=} \frac{\sqrt{N^{\alpha\beta}}}{\sqrt{\log\log N}} \left[ \alpha \eta^{\alpha - 1} \left( \frac{\xi(N)}{N^{\beta}} - \eta \right) + o\left( \left| \frac{\xi(N)}{N^{\beta}} - \eta \right| \right) \right]. 
\eal
Applying Proposition \ref{ClustersConcentration} and taking $\delta = \frac{1}{N^2}$, for every $i \in [d]$, there exists constants $c_i(\beta, \theta), C_i(\beta, \theta) > 0$ such that for $N > e^{C_i(\beta, \theta)}$,
\bal
P\left( \left| \frac{K_{i,N_i}(\beta, \theta)}{N_i^\beta} - S_{\beta, \theta, i} \right| \leq \frac{c_i(\beta, \theta)[\log\log(N_i + 2) + 2\log N]}{(N_i + \theta)^{\beta/2}} \right) \geq 1 - \frac{1}{N^2}. 
\eal
Since $\frac{N_i}{N} \to w_i$, there exists a constant $c_i'(\beta, \theta)$, which also depends on $w_i$, such that $N_i \geq c_i N$ for all $N$ sufficiently large. It follows that
\bal
P\left( \left| \frac{K_{i,N_i}(\beta, \theta)}{N_i^\beta} - S_{\beta, \theta, i} \right| > \frac{c_i''(\beta, \theta)\log N}{N^{\beta/2}} \right) \leq \frac{1}{N^2}
\eal
for some constant $c_i''(\beta, \theta)$ which depends on $\beta, \theta, w_i$ and for all $N$ sufficiently large. Now since
\bal
\sum_{N=1}^\infty P\left( \left| \frac{K_{i,N_i}(\beta, \theta)}{N_i^\beta} - S_{\beta, \theta, i} \right| > \frac{c_i''(\beta, \theta)\log N}{N^{\beta/2}} \right) \leq \sum_{N = 1}^\infty \frac{1}{N^2} < \infty, 
\eal
Borel-Cantelli implies that
\bal
\left| \frac{K_{i,N_i}(\beta, \theta)}{N_i^\beta} - S_{\beta, \theta, i} \right| \leq \frac{c_i''(\beta, \theta)\log N}{N^{\beta/2}}
\eal
almost surely for all $N$ sufficiently large. Furthermore, 
\bal
\sum_{N = 1}^\infty P\left( \bigcup_{i=1}^d \left\{ \left| \frac{K_{i,N_i}(\beta, \theta)}{N_i^\beta} - S_{\beta, \theta, i} \right| > \frac{c_i''(\beta, \theta)\log N}{N^{\beta/2}} \right\} \right) \leq \sum_{N = 1}^\infty \frac{d}{N^2} < \infty, 
\eal
so that by Borel-Cantelli again, we have that simultaneously for all $i \in [d]$, 
\bal
\left| \frac{K_{i,N_i}(\beta, \theta)}{N_i^\beta} - S_{\beta, \theta, i} \right| \leq \frac{c_i''(\beta, \theta)\log N}{N^{\beta/2}}
\eal
almost surely for all $N$ sufficiently large. In particular, we have that
\bal
\max_{i \in [d]}\left| \frac{K_{i,N_i}(\beta, \theta)}{N_i^\beta} - S_{\beta, \theta, i} \right| = O\left( \frac{\log N}{N^{\beta/2}} \right). 
\eal
Using this and the fact that $\left(\frac{N_i}{N}\right)^\beta \xrightarrow{a.s.} w_i^\beta$, we get
\bal
\left| \frac{\xi(N)}{N^\beta} - \eta \right| &= \left| \sum_{i=1}^d \left[ \left( \frac{N_i}{N} \right)^\beta \frac{K_{i,N_i}(\beta,\theta)}{N_i^\beta} - w_i^\beta S_{\beta, \theta, i} \right] \right| \\
&= \left| \sum_{i=1}^d \left[ \left( \frac{N_i}{N} \right)^\beta \left(\frac{K_{i,N_i}(\beta,\theta)}{N_i^\beta} - S_{\beta, \theta, i} \right) + \left( \left( \frac{N_i}{N} \right)^\beta - w_i^\beta \right) S_{\beta, \theta, i} \right] \right| \\
&\leq \sum_{i=1}^d  \left[ \left( \frac{N_i}{N} \right)^\beta \left| \frac{K_{i,N_i}(\beta,\theta)}{N_i^\beta} - S_{\beta, \theta, i} \right| + \left| \left( \frac{N_i}{N} \right)^\beta - w_i^\beta \right| S_{\beta, \theta, i} \right] \\
&= O\left( \frac{\log N}{N^{\beta/2}} \right).
\eal
Thus we get the upper bound 
\bal
&\left| \frac{\sqrt{N^{\alpha\beta}}}{\sqrt{\log\log N}} \left[ \alpha \eta^{\alpha - 1} \left( \frac{\xi(N)}{N^{\beta}} - \eta \right) + o\left( \left| \frac{\xi(N)}{N^{\beta}} - \eta \right| \right) \right] \right| \\
&\qquad \leq \frac{\sqrt{N^{\alpha\beta}}}{\sqrt{\log\log N}} \left| \frac{\xi(N)}{N^{\beta}} - \eta \right| \left( \alpha \eta^{\alpha - 1} + \frac{o\left( \left| \frac{\xi(N)}{N^{\beta}} - \eta \right| \right)}{\left| \frac{\xi(N)}{N^{\beta}} - \eta \right|} \right) \\
&\qquad = O\left( \frac{\sqrt{N^{\alpha\beta}}}{\sqrt{\log\log N}}  \frac{\log N}{N^{\beta/2}} \right) \\
&\qquad = O\left( \frac{\log N}{\sqrt{\log\log N}} \frac{1}{N^{\beta(1 - \alpha)/2}} \right),
\eal
so that
\bal
\left| \frac{\sqrt{N^{\alpha\beta}}}{\sqrt{\log\log N}} \left[ \alpha \eta^{\alpha - 1} \left( \frac{\xi(N)}{N^{\beta}} - \eta \right) + o\left( \left| \frac{\xi(N)}{N^{\beta}} - \eta \right| \right) \right] \right| \xrightarrow{a.s.} 0
\eal
as $N \to \infty$. 

Therefore we have shown that
\bal
\left| \frac{\sqrt{N^{\alpha\beta}} \left( \frac{\xi(N)}{N^{\beta}} - \eta \right)}{\sqrt{\log\log N}} \right| \xrightarrow{a.s.} 0
\eal
as $N \to \infty$, from which it follows that 
\bal
\frac{S_{\alpha, \theta_0}}{\sqrt{2}}\frac{\sqrt{N^{\alpha\beta}} \left( \left(\frac{\xi(N)}{N^{\beta}}\right)^\alpha - \eta^\alpha \right)}{\sqrt{\log\log N}} \xrightarrow{a.s.} 0
\eal
as $N \to \infty$. Combining the above finishes the proof of Case (4). 
\end{proof}


\section{Discussion and Final Remarks} \label{Sec:FinalRemarks}

We investigated the large $N$ asymptotic behavior of the number ${\bf K}_N$ of distinct values and the number ${\bf \widetilde M}_{\ell,N}$ of distinct values represented by $1 \leq \ell \leq N$ first-level clusters in a random sample of size $N$ from a hierarchical Pitman-Yor process. Our results show the impact of different levels in the hierarchy through the scalings of Gaussian fluctuations and the rate functions in the large deviation principles. Similar results are expected for more general hierarchical species sampling models, assuming the random index representation holds.  While we focused on hierarchical models involving two levels, we expect similar results for more levels.

\subsection{Further Asymptotic and Non-asymptotic Analysis of ${\bf K}_N$ and ${\bf \widetilde M}_{\ell,N}$}

Our study leaves open the problem of establishing a suitable Gaussian fluctuation for ${\bf \widetilde M}_{\ell,N}$. In particular, assuming the setting of Theorem \ref{HEPCLT}, we conjecture that as $N \rightarrow \infty$
\begin{equation}\label{gaus_fluc_r}
\sqrt{N^{\alpha\beta}}\left(\frac{{\bf \widetilde M}_{\ell,N}(\alpha,\theta_{0},\beta,\theta)}{N^{\alpha\beta}}-p_{\alpha}(r)S_{\alpha,\theta_{0}}\eta^{\alpha}\right)\stackrel{D}{\longrightarrow}\sqrt{p_{\alpha}(r)S_{\alpha,\theta_{0}}\eta^{\alpha}}N(0,1),
\end{equation}
where $p_{\alpha}(r)=\frac{\alpha(1-\alpha)^{(r-1)}}{r!}$. This remains an open problem even in the non-hierarchical Ewens-Pitman species sampling setting, where a preliminary result is given in \cite[Theorem 3.3]{BF24}. Establishing quantitative versions of Theorem \ref{HEPCLT}, as well as of \eqref{gaus_fluc_r} is also an interesting direction for future work. More generally, proving non-asymptotic results for the number of clusters and the number of clusters with frequency $r$, such as Berry-Esseen type inequalities and concentration inequalities, remains an open problem for both hierarchical and non-hierarchical species sampling models. We refer to \cite{DF20,DF21} and \cite{Per(22)} for some results under the Ewens-Pitman sampling model with $\alpha\in[0,1)$ and $\theta>-\alpha$.

Another research direction is related to the problem of estimating the parameters of hierarchical Ewens--Pitman sampling model, namely $\alpha, \beta \in (0,1)$, $\theta_0 > -\alpha$, and $\theta > -\beta$. In particular, the estimation of $\alpha$ and $\theta_{0}$ is known to be a challenging task  \cite{Bac(22)}. From Remark \ref{estim_alpha},
\begin{equation}\label{estimator_alfa}
\frac{{\bf \widetilde M}_{1,N}(\alpha, \theta_0, \beta, \theta)}{{\bf K}_N(\alpha, \theta_0, \beta, \theta)}
\end{equation}
is a (strongly) consistent estimator of $\alpha\in(0,1)$. An analogous result holds for the Ewens-Pitman sampling model \cite[Lemma 3.11]{Pit06}. One may further investigate \eqref{estimator_alfa}, as well as its non-hierarchical counterpart, by establishing a central limit theorem and, consequently, deriving a confidence interval for $\alpha$. To date, a consistent estimator of $\theta_{0}$ is unknown \cite{Bal(24)}.

It is worth mentioning a further line of research, recently initiated in the non-hierarchical setting of the Ewens-Pitman sampling model, concerning the large $n$ asymptotic behavior of $K_{n}(\alpha,\theta)$ and $M_{r,n}(\alpha,\theta)$ when the parameter $\theta$ grows with $n$, e.g. $\theta=\lambda n^{\beta}$, for some $\lambda>0$ and $\beta\geq0$ \cite{Cont25_1,Cont25_2}; see also \cite{Fen07} for early results on large deviations for $K_{n}$ in the Ewens sampling model. This direction has interesting applications in the context of microclustering, where the number of clusters grows with the sample size while their average size remains small, a feature commonly observed in modern large-scale data settings such as entity resolution and network analysis \cite{Ber25}. It would be natural to address the same question in the context of the hierarchical Ewens-Pitman sampling model.

\subsection{Bayesian Nonparametric Inference for Species Sampling Problems}

Proposition \ref{HEPLLN} and Proposition \ref{HEPblockLLN} pave the way for future research directions in Bayesian nonparametric statistics for {\em species sampling problems} \cite{Bal(24)}. In this context, one of the major problems is the estimation (or prediction) of the number of unseen species in sample. Specifically, given $N\geq1$ observed samples $\{(X_{i,1},\ldots, X_{i,N_i})\}_{i \in [d]}$ from a hierarchical Pitman-Yor process, the so-called unseen-species problem calls for estimating
\begin{displaymath}
{\bf{K}}_{M}^{(N)}=|\{(X_{i,N_{i}+1},\ldots, X_{i,N_i+M_{i}})\}_{i \in [d]}\setminus \{(X_{i,1},\ldots, X_{i,N_i})\}_{i \in [d]}|,
\end{displaymath}
namely the number of hitherto unseen species (symbols) that would be observed in $M=\sum_{i=1}^{d}M_{i}$ additional samples $\{(X_{i,N_{i}+1},\ldots, X_{i,N_i+M_{i}})\}_{i \in [d]}$ from the same model. Then, the interest is in the conditional distribution of ${\bf{K}}_{M}^{(N)}$ given $\{(X_{i,1},\ldots, X_{i,N_i})\}_{i \in [d]}$, i.e. the posterior distribution, whose expectation provides a Bayesian estimator of ${\bf{K}}_{M}^{(N)}$. Differently from the unseen-species problem in the non-hierarchical setting \cite{Lij(07),FLMP09}, obtaining closed-form expressions for the estimator of ${\bf{K}}_{M}^{(N)}$ is not possible, due to the analytical complexity inherent to the hierarchical construction \cite[Section 6]{CLOP19}. Instead, one may consider the problem of obtaining a large $M$ approximate estimator by characterizing the almost sure or $L^{p}$ large $M$ asymptotic behavior of ${\bf{K}}_{M}^{(N)}$, in analogy to Proposition \ref{HEPLLN}. A similar question holds with respect to the estimation of the number of hitherto unseen species with frequency $r$ that would be observed in $M=\sum_{i=1}^{d}M_{i}$ additional samples.


\section*{Acknowledgements}

The authors are grateful to two anonymous Referees for their careful reading of the manuscript, and for their helpful comments that have improved the paper.


\Address

\end{document}